\documentclass{amsart}

\usepackage{graphicx}%
\usepackage{multirow}%
\usepackage{amsmath,amssymb,amsfonts}%
\usepackage{amsthm}%
\usepackage{mathrsfs}%
\usepackage[title]{appendix}%
\usepackage{xcolor}%
\usepackage{textcomp}%
\usepackage{manyfoot}%
\usepackage{booktabs}%
\usepackage{algorithm}%
\usepackage{algorithmicx}%
\usepackage{algpseudocode}%
\usepackage{listings}%
\usepackage{amssymb}  
\usepackage{amsthm} 
\usepackage{amsmath}
\usepackage{amsfonts}
\usepackage{graphicx}
\usepackage{placeins}
\usepackage[utf8]{inputenc}
\usepackage[english]{babel}
\usepackage[all]{xy}

\usepackage{amsfonts}
\usepackage{amscd}
\usepackage{amsmath,amssymb,amsthm}
\usepackage{amstext}
\usepackage{amscd}
\usepackage{multicol}
\usepackage{graphicx}
\usepackage[all]{xy}
\usepackage{pstricks,pst-node,pst-tree}
\usepackage{dsfont}
\usepackage[pagewise]{lineno}

\newtheorem{lem}{Lemma}[section]
\newtheorem*{main1}{Theorem 1}
\newtheorem*{main2}{Theorem 2}
\newtheorem*{main3}{Theorem 3}

\newtheorem{coro}[lem]{Corollary}
\newtheorem{proposition}[lem]{Proposition}
\theoremstyle{definition}

\newtheorem{definition}[lem]{Definition}

\newtheorem{rem}[lem]{Remark}

\newtheorem{exs}[lem]{Examples}

\newcommand{\diam}{\mathop{\mathrm{diameter}}}

\newcommand{\C}{\mathbb{C}}
\newcommand{\R}{\mathbb{R}}

\newcommand{\N}{\mathbb{N}}
\newcommand{\Z}{\mathbb{Z}}

\newcommand{\bD}{\mathbb{D}}

\newcommand{\cK}{\mathcal{K}}

\newcommand{\bS}{\mathbb{S}}

\newcommand{\bT}{\mathbb{T}}
\newcommand{\bB}{\mathbf{B}}
\newcommand{\bC}{\mathbf{C}}
\newcommand{\bO}{\mathbf{O}}
\newcommand{\ie}{\emph{i.e},\,\,}

\DeclareMathOperator{\Int}{Int}

\author{Gabriela Hinojosa, Alberto Verjovsky and Juan Pablo D\'iaz}

\address{Centro de Investigaci\'on en Ciencias. Instituto de Investigaci\'on en Ciencias B\'asicas y Aplicadas. Universidad Aut\'onoma del Estado de Morelos. Av. Universidad 1001, Col. Chamilpa.
Cuernavaca, Morelos, M\'exico, 62209. }
\email{gabriela@uaem.mx}

\address{Instituto de Matem\'aticas Unidad
 Cuernavaca, Av. Universidad s/n. Col. Lomas de Chamilpa C\'odigo Postal 
 62210, Cuernavaca, Morelos.}
 
\email{albertoverjovsky@gmail.com}

\address{Centro de Investigaci\'on en Ciencias. Instituto de Investigaci\'on en Ciencias B\'asicas y Aplicadas. Universidad Aut\'onoma del Estado de Morelos. Av. Universidad 1001, Col. Chamilpa.
Cuernavaca, Morelos, M\'exico, 62209. }
\email{juanpablo.diaz@uaem.mx}

\thanks{This work was partially supported by PAPIIT (Universidad 
Nacional Aut\'onoma de M\'exico) project \#IN103324.}
\title[Beaded necklaces]{$N$-dimensional beaded necklaces and higher dimensional wild knots, invariant by a Schottky group} 

\begin{document} 
\maketitle 
\rightline{\footnotesize Dedicated to Francisco Javier Gonz\'alez Acu\~na (Fico)}
\rightline{\footnotesize  and Jos\'e Mar\'ia Montesinos Amilibia.}

\begin{abstract} Starting with a smooth, non-trivial $n$-dimensional knot $K\subset\bS^{n+2}$, and a beaded $n$-dimensional necklace subordinated to $K$, we construct a wild knot with a Cantor set
of wild points (\ie the knot is not locally flat in these points). The construction uses the conformal Schottky group acting on $\bS^{n+2}$, generated by inversions on the spheres which are the boundary of the ``beads''. We show that if $K$ is a fibered knot, then the wild knot is also fibered. We also study cyclic branched coverings along the
wild knots. This work generalizes the result presented in \cite{DH}.\\
\end{abstract}



\section{Introduction} \label{sec1}

The Cantor set was discovered by Henry John Stephen Smith in 1874 and mentioned by  Georg Cantor in 1883; the Sierpiński carpet and triangle are plane fractals first described by Wacław Sierpi\'nski in 1916, both generalizations of the Cantor set to two dimensions. All of the former examples helped lay the foundations of modern point-set topology. They led to notions such as local connectedness, different types of homology, Hausdorff dimension, and Fractal sets, among others. In the 1920s, with the works of Alexander,
Antoine, Artin, and Fox, among others, have described new examples of ``wild topology", such as wild knots and ``knotted" Cantor sets. In this century, this subject has seen great progress, in particular, thanks to the
use of computers and highly developed computer graphics; the subject has been popularized, and even the images created are part of modern art. Other examples of achievements in this century are, for instance, the result by Montesinos (\cite{Montesinos}, \cite{Montesinos1}) which states that every closed orientable 3-manifold is a 3-fold branched covering of $\mathbb{S}^{3}$ with branched set a wild knot.  Montesinos also  showed (\cite{M}) that  there exist uncountably many inequivalent wild knots whose cyclic branched coverings are $\mathbb{S}^3$,  via constructing wild disks in $\mathbb{S}^3$ such that each boundary has a tame Cantor set of locally wild points, and the corresponding $n$-fold cyclic covers of 
$\mathbb{S}^3$ branched over the boundary of each disk, are all $\mathbb{S}^3$. \\

\noindent In \cite{DH}, the first and second named authors extended the construction of  cyclic branched covers of the 3-sphere along  
{\it wild knots of dynamically defined type}; \ie 1-knots obtained via the action of a Kleinian group, in such a way that they construct a sequence of nested pearl chain necklaces, whose inverse limit space, that is the intersection space of this sequence of spaces, is a wild knot of dynamically defined type (compare Section 2).  Even more, they exhibited an example of a wild knot of dynamically defined type whose  $n$-fold cyclic branched cover of $\mathbb{S}^3$ along it is $\mathbb{S}^3$, for every $n\geq 2$. However, these wild knots are not the limit set of the corresponding Kleinian group (compare dynamically defined wild knots in \cite{Hin}); moreover, the corresponding Kleinian limit sets are Cantor sets, which are their set of wild points (compare Lemmas \ref{knot} and \ref{wildpoints}). We would like to recall that if $\Gamma$ denotes a Kleinian group, then its limit set consists of all the points of the 3-sphere that are accumulation points of some orbit of $\Gamma$ (for more details see \cite{maskit}). \\

\noindent In this work, we will generalize the construction of wild knots of dynamically defined type to higher dimensions, and we will study cyclic branched covers of the $(n+2)$-sphere along {\it wild $n$-knots invariant by a Schottky group}. 
An example of a wild 2-sphere in $\mathbb{S}^4$ which is the limit set of a geometrically
finite Kleinian group was obtained by the second-named author \cite{hinojosa2} and, independently, by Belegradek \cite{belegradek} (see also \cite{apanasov} for a wild limit set $\mathbb{S}^{2}\rightarrow \mathbb{S}^{3}$). In 2009, M. Boege, G. Hinojosa, and A. Verjovsky \cite{BHV}
obtained examples of wild $n$-spheres in $\mathbb{S}^{n+2}$ ($n=1,2,3,4$ and $5$) which are the limit set of geometrically
finite Kleinian groups. However, our desired wild $n$-knots are not the limit set of a Kleinian group; they are going to be the limit space of a nested sequence of beaded necklaces, which are obtained by the action of a Schottky group, and again, they will have a Cantor set of wild points. 
Next, we will give a brief description of them.\\

\noindent  Given a point $x\in\mathbb{S}^m$, let $B^m_r(x)=\lbrace{z}\in\bS^{m}: d(x,z)\leq{r}\rbrace$ denote the round, closed $m$-ball of radius $r>0$ centered at $x$, where $d$ is the standard metric induced from $\mathbb{R}^{m+1}$. Let $K\subset \mathbb{S}^{n+2}$ be a smooth $n$-dimensional knot. For $k\geq3$, consider $k$ points  $x_j\in{K}$, $j\in\lbrace1,\dots,k\rbrace$ with $r_j>0$ sufficiently small so that the balls $B^{n+2}_j:=B^{n+2}_{r_j}(x_j)$ are all disjoint and each pair
 $(B^{n+2}_j,B^{n+2}_j\cap{K})$ is a trivial $n$-tangle. An {\it $(n+2)$-dimensional $k$-beaded necklace $B(K,\bB_0)$ subordinated to $K$}, is a subset of $\bS^{n+2}$ defined as  $B(K,\bB_0)=K\cup{\bB_0}$, where 
$\bB_0=\overset{k}{\underset{j=1}\bigcup}\,B^{n+2}_j$ (see Definition \ref{bn}). \\

\noindent Let $\Gamma_{\bB_0}$ be the group generated by inversions $I_{j}$ 
through the $(n+1)$-sphere
$\Sigma^{n+1}_{j}:=\partial B^{n+2}_j$. We are interested in obtaining an $n$-dimensional wild knot invariant by a Schottky group, so that we will construct a nested sequence of beaded necklaces via the action of $\Gamma_{\bB_0}$, such that the corresponding inverse limit space
${\mathcal{K}}$  will be our desired wild $n$-knot. Even more, since $k\geq 3$, then ${\mathcal{K}}$ will have a Cantor set of wild points (for more details see Section \ref{dyn}). \\

\noindent In this paper, we prove some topological properties of higher-dimensional wild knots invariant by a Schottky group. Among them stand out the following.

\begin{main1}
Let $B(K,\bB_0)$ be a $k$-beaded necklace subordinated
to the non-trivial, smooth fibered $n$-knot $K$. Let $\Gamma_{\bB_0}$ be the group
generated by inversions through each $(n+1)$-sphere $\Sigma^{n+1}_{j}=\partial B^{n+2}_{j}$, where $B^{n+2}_{j}\subset \bB_0$, $j=1,2,\ldots,k$, and consider the 
inverse limit space  $\mathcal{K}$.  Then:
\begin{enumerate}
\item There exists a locally trivial fibration $\psi
:(\mathbb{S}^{n+2}\smallsetminus\mathcal{K})\rightarrow\mathbb{S}^{1}$, where the
 fiber $\Sigma_{\theta}=\psi^{-1}(\theta)$ is an orientable $(n+1)$-manifold with one end,
and if the fiber of  $K$ has non-trivial homology in dimension $r$,
then $H_r(\psi^{-1}(\theta),\Z)$ is infinitely generated.
\item  $\overline{\Sigma_{\theta}}{\smallsetminus}\Sigma_{\theta}=\mathcal{K}$, where
$\overline{\Sigma_{\theta}}$ is the closure of $\Sigma_{\theta}$ in $\mathbb{S}^{n+2}$.
\end{enumerate}
\end{main1}

\noindent Consider two $k$-beaded necklaces  $B(K,\bB_0)$ and $B(L,\bC_0)$ subordinated to the $n$-knots $K$ and $L$ respectively, such that $\bB_0=\bigcup_{j=1}^n B_j$ and $\bC_0=\bigcup_{i=1}^n C_i$ where $B_j$, $C_i$ are closed $(n+2)$- balls $(i,\,j=1,2,\dots,k)$. If they are equivalent (see Definition \ref{NE}), we have the following.

\begin{main2}\label{main2}
Let $B(K,\bB_0)$ and $B(L,\bC_0)$ be two equivalent $k$-beaded necklaces subordinated to the smooth $n$-knots $K$ and $L$, respectively. Then the corresponding inverse limit spaces $\mathcal{K}$ and $\mathcal{L}$ are equivalent wild $n$-knots.
\end{main2}

\noindent Moreover, we extend the construction of cyclic branched covers to this kind of wild knots, and we show the following.

\begin{main3}\label{main3}
Let $B(K,\bB_0)$ be a $k$-beaded necklace subordinated
to the non-trivial, smooth $n$-knot $K$, and let  $\mathcal{K}$ be the corresponding wild $n$-knot. Then, for each integer $q$, there exists a $q$-fold cyclic branched cover $\Psi:\mathbb{M}_q\rightarrow \mathbb{S}^{n+2}$ along $\mathcal{K}$ such that $\mathbb{M}_q$ is a compact and connected space.
\end{main3}

\noindent We would like to remark that it is a very difficult problem to prove that the Freudenthal compactification of a non-compact $(n+2)$-manifold, whose space of 
Freudenthal ends is a Cantor set, is again an $(n+2)$-manifold, or even a homotopy manifold (for more details see \cite{fox} and \cite{M}). This is not always the case.

\section{Construction of $n$-dimensional wild knots invariant by a Schottky group}\label{dyn}

In this section, we will describe the construction of $n$-dimensional wild knots invariant by a Schottky group generated by inversions on disjoint spheres. These are obtained as the inverse limit of a nested sequence of beaded necklaces. We will start with some previous definitions.\\

\noindent In classical knot theory, a subset $K$ of a space $X$ is a {\it knot} if $K$ is homeomorphic to a sphere
$\mathbb{S}^{p}$. Two knots $K$, $K'$ are {\it equivalent} if there is a homeomorphism $h:X\rightarrow X$ such that $h(K)=K'$;
in other words $(X,K)\cong (X,K')$. However, a knot $K$ is sometimes defined to be an embedding
$K:\mathbb{S}^{p}\rightarrow\mathbb{S}^{n}$ (see \cite{mazur}, \cite{rolfsen}).
We shall also find this convenient at times and will use the same symbol to
denote either the map $K$ or its image $K(\mathbb{S}^{p})$ in $\mathbb{S}^{n}$.\\

\noindent For an integer $m$, let $\bS^m$ denote the unit $m$-sphere in $\R^{m+1}$ with its standard metric $d$. For $x\in\mathbb{S}^m$, 
let $B^m_r(x)=\lbrace{z}\in\bS^{n}: d(x,z)\leq{r}\rbrace$ be the round, closed $m$-ball of radius $r>0$ centered at $x$. In particular, if $x$ is the origin and $r=1$,
we will have the  $m$-disk $\mathbb{D}^{m}:=B^m_1(0)$.

\begin{definition}
 We say that a point $x\in K$ is 
{\it locally flat} or {\it locally tame} 
if there exists an open neighborhood $U$ of $x$ such that
there is a homeomorphism of pairs: 
$(U,U\cap K)\sim (\Int(\mathbb{D}^{n+2}),\Int(\mathbb{D}^{n}))$.
Otherwise, $x$ is said to be a {\it wild} point of $K$. 
An $n$-knot $K$ is {\it locally flat} or {\it locally tame} if all of its
points are locally flat. Otherwise, we say $K$ is a {\it wild knot}.
\end{definition}

\begin{definition}\label{ntangle}
An {\it oriented n-dimensional tame single-strand tangle} is a couple $D=({B}^{n+2},T)$ satisfying the following conditions:
\begin{enumerate}
\item ${B}^{n+2}$ is homeomorphic to the $(n+2)$-disk $\mathbb{D}^{n+2}$, and $T$ is homeomorphic to the $n$-disk $\mathbb{D}^{n}$.
\item The pair $({B}^{n+2},T)$ is a proper manifold pair; \ie $\partial T\subset \partial {B}^{n+2}$ and
$\mbox{Int}(T)\subset \mbox{Int}({B}^{n+2})$.
\item $(B^{n+2},T)$ is locally flat (see \cite{rushing}, p.33, for the notion of local flatness for pairs $(M,N)$ of manifolds with boundary  with $N\subset{M}$ and 
$\partial{N}\subset\partial{M}$).
\item ${B}^{n+2}$ has an orientation which induces the canonical orientation on its boundary
$\partial B^{n+2}$.
\item $(\partial B^{n+2},\partial T)$ is homeomorphic to $(\partial\mathbb{D}^{n+2},\partial\mathbb{D}^{n})=(\bS^{n+1},\bS^{n-1})$.
\end{enumerate}

\noindent {\bf Meridian.}\label{meridian} Write $X_T := B^{n+2} \smallsetminus \nu(T)$ for the complement of an open tubular neighbourhood $\nu(T)\cong T\times \mathring{\mathbb{D}}^2$, where 
$\mathring{\mathbb{D}}^2$ denotes the interior of $\mathbb{D}^2$. A \emph{meridian} of a component of $T$ is a loop generating the $S^1$ of the $\partial \nu(T)\cong T\times \bS^1$ factor.
\end{definition}

\noindent Compare to Zeeman's definition of ball-pair in (\cite{zeeman}).\\

\noindent Two oriented tangles $D_{1}=({B}_{1}^{n+2},T_{1})$, $D_{2}=({B}_{2}^{n+2},T_{2})$ are \emph{equivalent} if there exists an orientation-preserving
homeomorphism of ${B}_{1}^{n+2}$ onto ${B}_{2}^{n+2}$ that sends $T_{1}$ to $T_{2}$. A tangle is \emph{unknotted} if it is
equivalent to the trivial tangle $(\mathbb{D}^{n+2},\mathbb{D}^{n})$ (for more details see \cite{BHV}).\\

\noindent Given an oriented tangle $D=({B}^{n+2},T)$, the pair $(\partial B^{n+2},\partial T)$ is homeomorphic to the pair
$(\mathbb{S}^{n+1},\mathbb{S}^{n-1})$, via a homeomorphism $f$.
Then $D$ determines canonically a knot $K\subset\mathbb{S}^{n+2}$, in the following way:
$(\mathbb{S}^{n+2},K)=({B}^{n+2},T)\cup_{f}(\mathbb{D}^{n+2},\mathbb{D}^{n})$. Conversely, given a smooth knot $K\subset\mathbb{S}^{n+2}$, there exists a smooth ball $B^{n+2}$ such that $(B^{n+2},B^{n+2}\cap K)$ is equivalent to the trivial tangle. The tangle
\begin{equation}\label{canonicaltangle}
K_{T}=(\mathbb{S}^{n+2}\smallsetminus\mbox{Int}(B^{n+2}),K\smallsetminus\mbox{Int}(B^{n+2}\cap K))
\end{equation}
\noindent is called the {\it
canonical tangle} associated to $K$. Notice that if $K$ is not the trivial knot, then $K_{T}$ is not equivalent to
the trivial tangle. In this case, we say that $K_{T}$ is knotted.

\begin{definition}\label{connectedsum} 
The {\it connected sum} of the oriented tangles $D_{1}=({B}_{1}^{n+2},T_{1})$ and $D_{2}=({B}_{2}^{n+2},T_{2})$ for $n>1$,
denoted by $D_{1}\# D_{2}$,
can be defined as follows:  Since  $D_{i}$ is locally flat $i=1,2$, there exist sets $U_{i}\subset \partial B_{i}$ closed in $\partial B_i$, 
such that $\mbox{Int}(U_{i})\cap \partial T_{i}\neq\emptyset$ and, the pair $(U_{i}, U_{i}\cap T_{i})$ is homeomorphic to
$(\mathbb{D}^{n+1},\mathbb{D}^{n-1})$. Choose an orientation-reversing homeomorphism of pairs
$(U_{1}, U_{1}\cap T_{1})$ onto
$(U_{2}, U_{2}\cap T_{2})$, (\ie $h:U_{1}\to U_{2}$ is a homeomorphism
and the restriction, $g$, of $h$ to $U_{2}\cap T_{2}$, is a homeomorphism from
$U_{1}\cap T_{1}$ to $U_{2}\cap T_{2}$).
Then, the connected sum $D_{1}\# D_{2}$ is defined as the pair:
\[
D_{1}\# D_{2}=({B}_{1}^{n+2},T_{1})\cup_{h} ({B}_{2}^{n+2},T_{2})
=({B}_{1}^{n+2}\cup_{h}{B}_{2}^{n+2}, T_1\cup_{g}T_2)\]
where $({B}_{1}^{n+2},T_{1})\cup_{h} ({B}_{2}^{n+2},T_{2})$ is the space obtained
from the disjoint union of ${B}_{1}^{n+2}$ and ${B}_{2}^{n+2}$
identifying $U_1$ with $U_2$ by the homeomorphism $h$ and $T_1\cup_{g}T_2$ is obtained from the disjoint union of $T_1$ and $T_2$ by identifying 
$U_{1}\cap T_{1}$ with $U_{2}\cap T_{2}$ via the homeomorphism $g$.
Then, ${B}_{1}^{n+2}\cup_{h}{B}_{2}^{n+2}:\overset{def}={B}_{1}^{n+2}\#{B}_{2}^{n+2}$ is homeomorphic to $\bD^{n+2}$ and
$T_1\cup_{g}T_2:\overset{def}=T_1\#T_2$ is homeomorphic to ${\bD}^{n}$ so that $D_{1}\# D_{2}$
is the tangle $({B}_{1}^{n+2}\#{B}_{2}^{n+2}, T_1\#T_2)$.
\end{definition}

\begin{rem}
The connected sum does not depend on the choice of the homeomorphism  $h$ and the sets $U_{i}$. 
The connected sum is an associative operation, so it makes sense to speak
of the connected sum of a finite number of $n$-dimensional tangles.
\end{rem}

\noindent One has the following proposition, which is an easy consequence of van Kampen's theorem:

\begin{proposition}[Fundamental group of a connected sum of tangles]\label{vankampen}
\label{thm:amalgamation}
Let $n\ge1$ and $r\geq2$. Let $D_1=(B_1^{n+2},T_1),\dots,(B_r^{n+2},T_r)$ be a collection of $r$ $n$-tangles as above. Let $D_1\#D_2\#\dots\#{D_r}$ be their connected sum.
Let $X_1=B_1^{n+2}\smallsetminus{T_1},\dots,X_r=B_r^{n+2}\smallsetminus{T_r}$ and $X := {B}_{1}^{n+2}\#\dots\#{B}_{r}^{n+2}\smallsetminus(T_1\#\dots\#T_r)$ 
denote the respective tangle complements. Then:
\[
\pi_1(X)\ \cong\ (\pi_1(X_1)*\dots *\, \pi_1(X_r))/F_r,\]
the free product of $\pi_1(X_1),\dots,\pi_1(X_r)$ of the $r$ fundamental groups of the complements over a free group $F_r$ of rank $r$. 
Under this isomorphism the canonical inclusions $F_r\hookrightarrow \pi_1(X_i)$ send a fixed free basis 
$\{\mu_1,\ldots,\mu_r\}$ to the classes of the $r$ chosen meridians in $\pi_1(X_i)$ encircling the local strand pieces used to perform the sum.

\noindent In particular, for $r=2$ one obtains an amalgam over $\mathbb{Z}$:
\[
\pi_1(X)\ \cong\ \pi_1(X_1)\ *_{\langle \mu_1=\mu_2\rangle}\ \pi_1(X_2),
\]
\ie the free product with the two chosen meridians identified.
\end{proposition}
%%%%

\begin{coro}
Let $K\subset \bS^{n+2}$ be an $n$-knot and let $\overline{K}$ denote its mirror image under an inversion
with respect to a round $(n+1)$-sphere in $\bS^{n+2}$.
Then the knot groups of the connected sums $K\#K$ and $K\#\overline{K}$ are isomorphic:
\[
\pi_1\big(\bS^{n+2}\setminus (K\#K)\big)\cong \pi_1\big(\bS^{n+2}\setminus (K\#\overline{K})\big).
\]
\end{coro}
%\vskip .2cm
\noindent We consider now the iterated connected sum of a knot with itself. Let \(K\subset \bS^{n+2}\) be a (tame) knot. Denote by \(G(K)=\pi_{1}\big(\bS^{n+2}\setminus K\big)\) the knot group of \(K\). For an integer \(r\ge 1\) write
\[
K^{\# r} := \underbrace{K\# K\#\cdots\# K}_{\text{$r$ copies}}
\]
for the connected sum of \(r\) copies of \(K\). Let \(m\in G(K)\) denote the homotopy class of a meridian loop of \(K\) (any choice of meridian gives a conjugate element; we fix one).

\begin{coro}\label{conn}
For every knot \(K\subset \bS^{n+2}\) and every integer \(r\ge 1\),
\begin{equation}\label{rconnected}
\pi_{1}\big(\bS^{n+2}\setminus K^{\# r}\big)
\;\cong\;
\underbrace{G(K)\;*_{\langle m\rangle}\;G(K)\;*_{\langle m\rangle}\;\cdots\;*_{\langle m\rangle}\;G(K)}_{\text{$r$ copies}},
\end{equation}
the iterated free product with amalgamation of \(r\) copies of \(G(K)\) where, in each amalgamation, the distinguished infinite cyclic subgroup \(\langle m\rangle\cong\mathbb{Z}\) generated by a meridian of that copy of \(K\) is identified with the same central amalgamating copy of \(\mathbb{Z}\).
\end{coro}

\noindent In words: the knot group of the connected sum is obtained by taking \(r\) copies of the knot group and identifying their meridian subgroups.

\begin{rem}\label{generatorsgrow} 
As a consequence of the previous results, we know that 
\[
\pi_{1}\big(\bS^{n+2}\setminus K^{\# r_1}\#(\overline{K})^{\# r_2}\big)\;\cong\;
\underbrace{G(K)\;*_{\langle m\rangle}\;G(K)\;*_{\langle m\rangle}\;\cdots\;*_{\langle m\rangle}\;G(K)}_{\text{$(r_1+r_2)$ copies}},
\]
It follows, in particular, that if $K\subset\bS^{n+2}$ is an $n$-knot with 
$\pi_1(K\smallsetminus\bS^{n+2})$ bigger than $\Z$ (\ie the group is not the infinite cyclic group represented by a meridian) then $\pi_{1}\big(\bS^{n+2}\setminus K^{\# r}\big)$ cannot have a set of generators of cardinality less than $r$.
\end{rem}

%%%%

\begin{figure}[ht]  
\begin{center}

%\hspace{.3cm}
\centering
\includegraphics[width=7cm, height=3.5cm ]{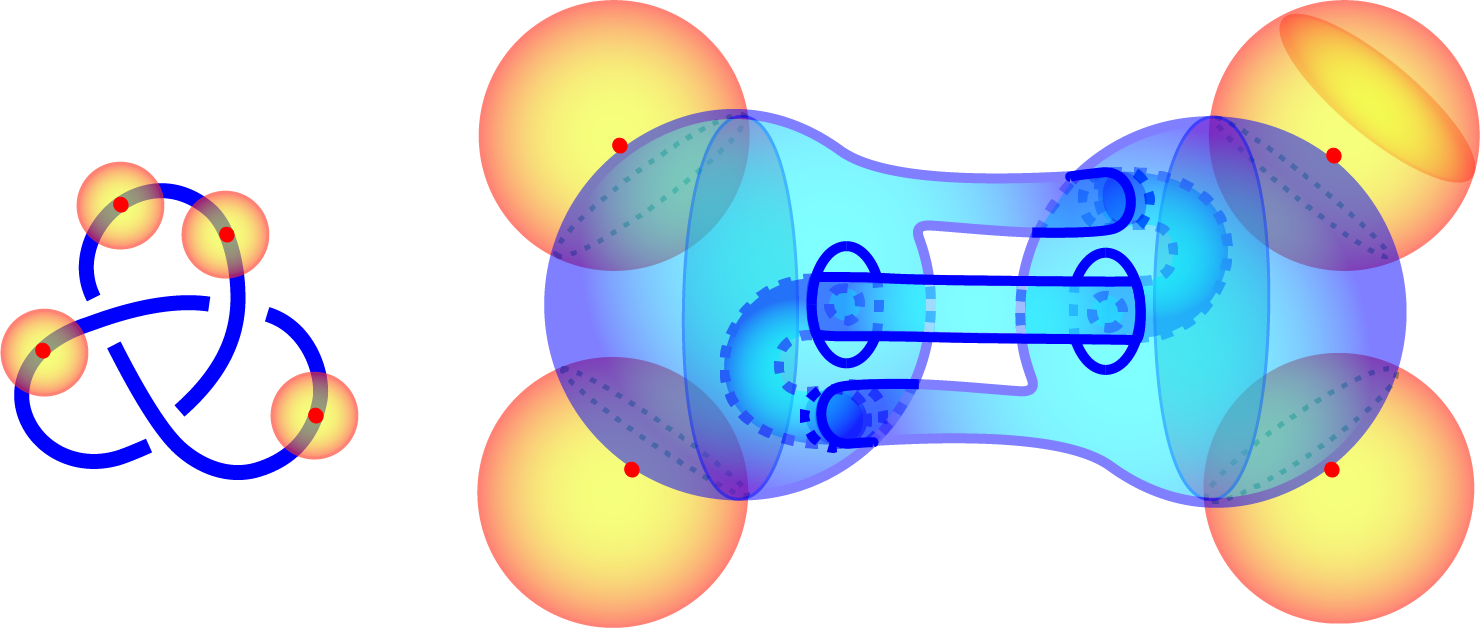}
\end{center}
\caption{\sl A beaded necklace of dimensions 1 and 2. In the schematic figure at the right, the blue part is the $n$-dimensional ``thread" and the yellow balls are the ``beads".}
\label{F1}
\end{figure}

\begin{definition}\label{bn}
Let $K\subset\bS^{n+2}$ be a smooth $n$-dimensional knot
\ie $K$ is the image of a smooth embedding $f:\bS^n\to\bS^{n+2}$. 
 An {\it $(n+2)$-dimensional beaded necklace $B(K,\bB_0)$ subordinated to $K$}, is a subset of $\bS^{n+2}$ which is a union $B(K,\bB_0)=K\cup{\bB_0}$, where 
$\bB_0=\overset{k}{\underset{j=1}\bigcup}\,B^{n+2}_{r_j}(x_j)$ is a finite union of disjoint, closed
$(n+2)$-balls $B^{n+2}_j:=B^{n+2}_{r_j}(x_j)$ called {\it beads} ($j\in\lbrace1,\dots,k\rbrace, k\geq3$), with the following
properties  (see Figure \ref{F1}):
\begin{enumerate}
\item $x_j\in{K}, \forall \, j\in\lbrace1,\dots,k\rbrace$

\item the radii $r_j$ are sufficiently small so that the pair 
\[(B^{n+2}_j,B^{n+2}_j\cap{K}),\]
is a trivial tangle of dimension $n$ (trivial $n$-tangle), for $j\in\lbrace1,\dots,k\rbrace$.
We refer to $K$ as the {\bf thread} of the beaded collar.
\end{enumerate}
\end{definition}

%%

\noindent We have the following lemma whose proof is elementary:

\begin{lem}[Beadded necklaces retract to their threads] \label{collars-retract}
Let $B(K,\bB_0)$ be a necklace subordinated to the $n$-dimensional knot $K$ (the ``thread"). Let $\bB_0=\overset{k}{\underset{j=1}\bigcup}\,B^{n+2}_{j}$ be the set of beads. Then $K$ is a strong deformation retract of $B(K,\bB_0)$. Furthermore, 
let $B^{n+2}_{i},\, i\in\lbrace1, \dots, k\rbrace$ be one of the beads, then 
${\hat B}^{n+2}_{i}\overset{def}=\bS^{n+2}\smallsetminus{\rm{Int}\,(B^{n+2}_{i}})$ is a round, closed, $(n+2)$-ball, and therefore
the pair $({\hat B}^{n+2}_{i}, {\hat B}^{n+2}_{i}\cap{K})$ is a tangle. In addition,
${\hat B}^{n+2}_{i}\cap{K}$ is a strong deformation retract of 
$\underset{{j\in\lbrace1, \dots, k\rbrace, \, j\neq{i}}}\bigcup\,{B^{n+2}_{j}}$.
\end{lem}

\noindent Let $\Gamma_{\bB_0}$ be the group generated by inversions $I_{j}$ through the $(n+1)$-sphere
$\Sigma^{n+1}_{j}:=\partial B^{n+2}_{j}$ ($j=1,\ldots,k$). Then $\Gamma_{\bB_0}$ is a discrete subgroup of
$\mbox{M\"ob}(\mathbb{S}^{n+2})$ whose limit set $\Lambda(\Gamma_{\bB_0})$ is a Cantor set, since $k\geq 3$ (\cite{kap1}). Even more clearly, its fundamental domain  is 
$D=\mathbb{S}^{n+2}\smallsetminus\bB_0$, hence
$\Gamma_{\bB_0}$ is a geometrically finite Kleinian group of Schottky type (\cite{kap1}, \cite{maskit}). \\

\noindent We are interested in constructing $n$-dimensional wild knots, so  we will construct a nested sequence of beaded necklaces via the action of $\Gamma_{\bB_0}$ into our beaded-necklace $B(K,\bB_0)$ to get an inverse limit space
${\mathcal{K}}$,  which will be our desired wild $n$-knot. \\ 

\begin{figure}[ht]  
\begin{center}
%\hspace{.3cm}
\includegraphics[height=3.5cm]{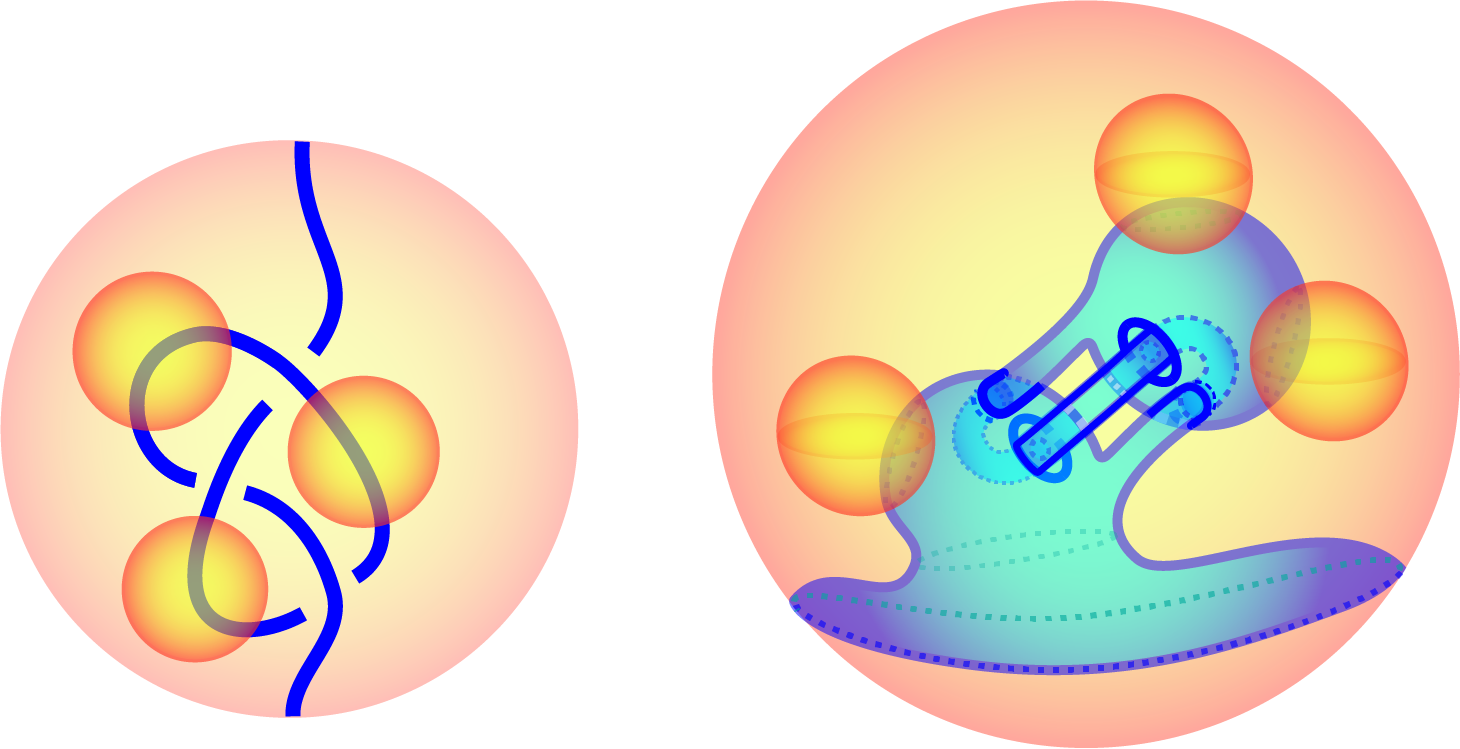}
\end{center}
\caption{\sl A schematic picture of the non-trivial $n$-tangle $(B^{n+2}_{j}, \kappa_{1_j})$ which contains a beaded  $n$-strand
$\kappa_{1_j}=I_j(B(K,\bB_0)\smallsetminus{B^{n+2}_{j}})$ consisting of the union of $k-1$ disjoint closed $(n+2)$-balls and an $n$-disk $\overline{K}'$.}
\label{strand}
\end{figure}

\noindent Notice that if we invert with respect to  $\Sigma^{n+1}_{j}$, both
an opposite oriented image (mirror image) of the closed $n$-disk $K\smallsetminus\Int (B^{n+2}_j\cap K)$, say $\overline{K}'$,  and the corresponding beaded  $(n+2)$-strand 
$B(K,\bB_0)\smallsetminus\Int (B^{n+2}_{j})$ are sent into the ball $B^{n+2}_{j}$. In other words,  $B^{n+2}_{j}$ contains a beaded  $(n+2)$-strand
$\kappa_{1_j}=I_j(B(K,\bB_0)\smallsetminus\Int (B^{n+2}_{j}))$, in such a way that if $K$ is a non-trivial knot, then $(B^{n+2}_{j}, \kappa_{1_j})$ is a non-trivial $n$-tangle. In fact, if we contract each ball of  $\kappa_{1_j}$ to a point, we get  $\overline{\kappa_{1_j}}$ such that the pair $(B^{n+2}_{j}, \overline{\kappa_{1_j}})$ is homeomorphic to the tangle
$(C,\overline{K}')$, where $C$ is a closed $(n+2)$-ball  (see Figure \ref{strand}). So, after this inversion, we obtain a new beaded necklace 
 \[
 B(K_{1_j},\bB_{1_j})=(B(K,\bB_0)\smallsetminus\Int(B^{n+2}_j))\cup_{I_j} \kappa_{1_j}, 
 \]which is gotten from $B(K,\bB_0)$ replacing the ball $B^{n+2}_j$ by $\kappa_{1_j}$, \noindent where  the 
 $n$-knot 
$K_{1_j}=(K\smallsetminus\Int(B^{n+2}_j\cap K))\cup_{I_j} I_j(K\smallsetminus\Int (B^{n+2}_j\cap K))$  is the connected sum of $K$ with its mirror image $\overline{K}$; \ie $K_{1_j}\cong K\#\overline{K}$.  Similarly, we obtain the corresponding set of  beads $\bB_{1_j}=(\bB_0\smallsetminus{B^{n+2}_j})\cup \kappa^{\circ}_{1_j}$, where 
$\kappa^{\circ}_{1_j}=I_j(\bB_0\smallsetminus{B^{n+2}_{j})}$.\\
 
\noindent Now, we invert with respect to $\Sigma^{n+1}_{r}$, $r\neq j$. As before, both a mirror copy of  the beaded  $(n+2)$-strand 
$B(K,\bB_0)\smallsetminus\Int(B^{n+2}_{r})$ and a mirror copy of the $n$-disk  
$K\smallsetminus\Int(B^{n+2}_r\cap K)$ go into $B^{n+2}_{r}$, obtaining a new beaded necklace 
\[
B(K_{1_{j,r}},\bB_{1_{j,r}})=(B(K_{1_j},\bB_{1_j})\smallsetminus\Int(B^{n+2}_r))\cup_{I_r} \kappa_{1_r},
\]
\noindent where $\kappa_{1_r}=I_r(B(K,\bB_0)\smallsetminus{B^{n+2}_{r}})$. This beaded necklace  is gotten from $B(K_{1_j},\bB_{1_j})$ replacing the ball 
$B^{n+2}_r$ by $\kappa_{1_r}$, and is subordinated to the $n$-knot 
\[
K_{1_{j,r}}=(K_{1_j}\smallsetminus\Int(B^{n+2}_r\cap K))\cup_{I_r} I_r(K_{1_j}\smallsetminus(B^{n+2}_r\cap K)),
 \]
 \noindent which is the connected sum of $K$ and two copies of its mirror image;
 \ie $K_{1_{j,r}}\cong K\#\overline{K}\#\overline{K}$.  Similarly, we obtain the corresponding set of  beads 
 \[
 \bB_{1_{j,r}}=(\bB_{1_j}\smallsetminus{B^{n+2}_r})\cup \kappa^{\circ}_{1_r},
 \]
\noindent where 
$\kappa^{\circ}_{1_r}=I_r(\bB_0\smallsetminus{B^{n+2}_{r}})$. \\

\begin{figure}[ht]  
\begin{center}
%\hspace{.3cm}
\includegraphics[height=3.5cm]{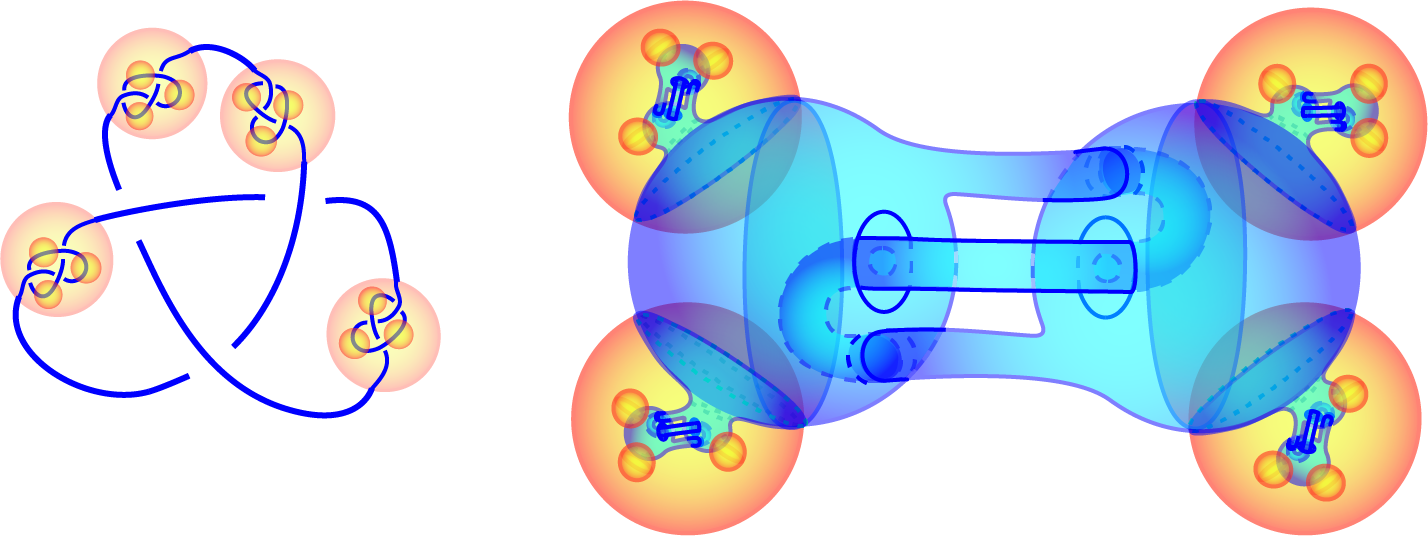}
\end{center}
\caption{\sl A schematic picture of a beaded necklace after the first stage of the inverting process. In the schematic figure at the right, the blue part is the $n$-dimensional ``thread" and the yellow spheres are the ``beads".} 
\label{F3}
\end{figure}

\noindent We continue in this way, after inverting with respect to each $\Sigma^{n+1}_{j}$ $j\in\{1,\dots, k\}$, we get
$\bB_{1}$ which consists of the union of $l_1=k(k-1)$ beads, $B^{1}_r$ $r\in\{1,\dots, l_1\}$, subordinated to a new $n$-knot $K_{1}$ which
is in turn isotopic to the connected sum of $K$ and $k$ copies of its mirror
image $\bar{K}$; \ie $K_{1}\cong K\#\overline{K}\#\overline{K}\#\ldots\#\overline{K}$ ($k$-times). More specific, for $j\in\{1,\ldots,k\}$, each $\kappa^{\circ }_{1_j}$ consists of the union of $k-1$ beads $B^1_{1_j}$ where the index $1_j$ belongs to a subset of the set of indexes $\{1,\ldots,k(k-1)\}$. Notice that for $l\in\{1,\ldots,k-1\}$ the balls $B^1_{(j-1)(k-1)+l}\subset B^0_j:=B_j$. Hence 
$\kappa^{\circ}_{1_j}=\bigcup_{l=1}^{k-1} B^{1}_{(j-1)(k-1)+l}$, so $\bB_{1}$ is the union of all $\kappa^{\circ}_{1_j}$. Observe that  all the $n$-knots determined from the tangles $(B^{n+2}_{j}, \overline{\kappa_{1_j}})$ are equivalent to the mirror image of $K$, hence $K_{1}\cong K\#\overline{K}\#\overline{K}\#\ldots\#\overline{K}$ ($k$-times). At the end of the first stage of the inverting process, we get a new beaded necklace  $B(K_1,\bB_1)=K_1\cup\bB_{1}$ such that $\bB_{1}\subset \bB_{0}$ and 
$B(K_1,\bB_1)\subset B(K,\bB_0)$ (see Figure~\ref{F3}). \\ 

\noindent We continue inductively,  so at the $(m-1)^{th}$ stage we have a set of beads $\bB_{m-1}$ which consists of $l_{m-1}=k(k-1)^{m-1}$ beads $B_j^{m-1}$, 
$j\in\{1,\dots, l_{m-1}\}$, and  is  subordinated to an 
$n$-knot $K_{m-1}$ which is in turn isotopic to the connected sum of $k^{m-2}+1$ addends, such that each addend is either the original $n$-knot $K$ or its mirror image 
$\overline{K}$ (recall that a composition of an even number or inversions is orientation-preserving). Let  $B(K_{m-1},\bB_{m-1})=K_{m-1}\cup\bB_{m-1}$ be the corresponding new beaded necklace, so $\bB_{m-1}\subset \bB_{m-2}$ and 
$B(K_{m-1},\bB_{m-1})\subset B(K_{m-2},\bB_{m-2})$. As in the first stage, we have that  each bead of $\bB_{m-2}$, 
say $B_j^{m-2}$ ($j\in\{1,\dots, l_{m-2}\}$ where $l_{m-2}=k(k-1)^{m-2}$), contains a  beaded strand
$\kappa_{(m-1)_j}$ consisting of the union of an $n-disk$ that determines an $n$-knot $\widehat{\kappa_{(m-1)_j}}$ equivalent to the original $n$-knot $K$ or  its mirror image $\overline{K}$ (depending if $m-1$ is even or odd), and a beaded strand $\kappa^{\circ}_{(m-1)_j}$ that is the union of $k-1$ disjoint beads that can enumerated using the
numeration (positions) coming from the previous level $k_{m-2}$ (see stage 1), namely
$B_{(j-1)(k-1)+1}^{m-1}, \dots , B_{j(k-1)}^{m-1}$ subordinated to the corresponding $n$-disk. Hence  the union of all the strands $\kappa^{\circ}_{(m-1)_j}$ is $\bB_{m-1}$ and the connected sum of $K_{m-2}$ with all the corresponding $n$-knots
 $\widehat{\kappa_{(m-1)_j}}$ is $K_{m-1}$.\\

\begin{figure}[ht] 
 \begin{center}
  \includegraphics[height=5cm]{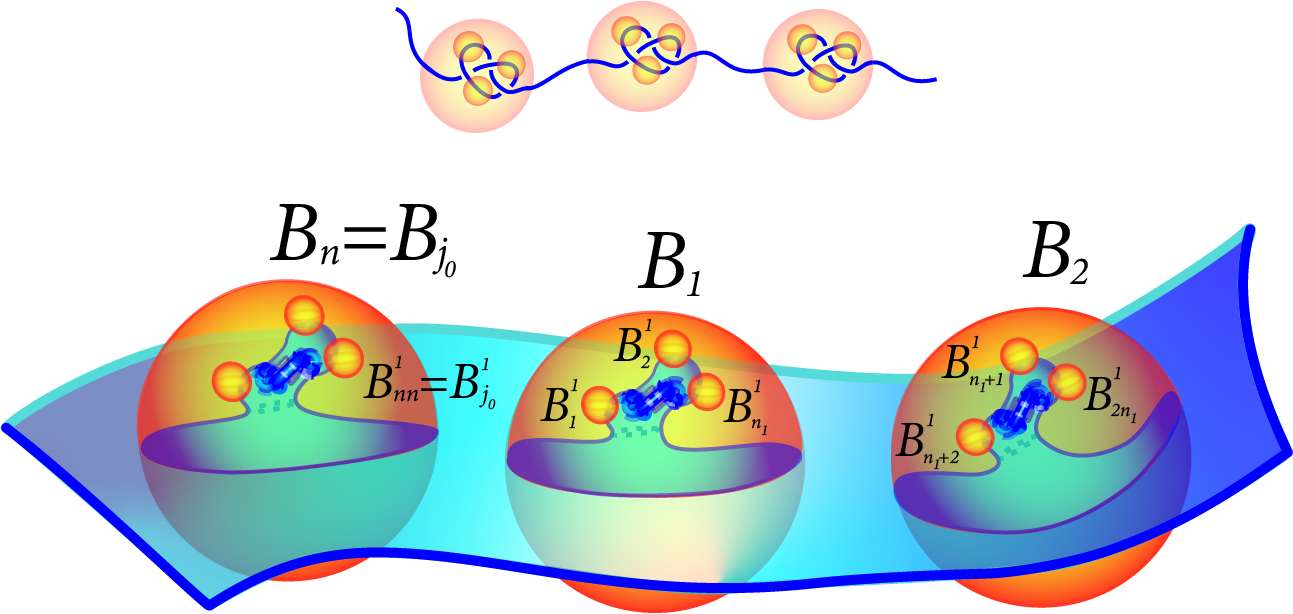}
\end{center}
\caption{\sl Schematic diagram of the nested beads. In the lower figure, the blue part is the $n$-dimensional ``thread" and the yellow spheres are the ``beads".} 
\label{Schem}
\end{figure}

\noindent The  $m^{th}$ stage of the inverting process is obtained from the previous stage after inverting with respect to each $\Sigma^{n+1}_{j}$. Then  we have a new  beaded set
$\bB_{m}$ which is the union of  $l_m=k(k-1)^{m}$ beads $B_j^{m}$, $j\in\{1,\dots, l_m\}$, subordinated to an
$n$-knot $K_{m}$. Let  $B(K_m,\bB_m)=K_m\cup\bB_{m}$ be the corresponding new beaded necklace, so by construction $\bB_{m}\subset \bB_{m-1}$ and 
$B(K_m,\bB_m)\subset B(K_{m-1},\bB_{m-1})$. Notice that the diameter of each bead $B_j^{m}$ tends to zero as $m$ goes to infinity (see Figure \ref{Schem}). \\

\noindent Therefore,  the inverse limit  space $\mathcal{K}$ is  given by  
\begin{equation}\label{beadsintersection}
\mathcal{K}=\varprojlim_{m} B(K_m,\bB_m)=\bigcap_{m=0}^{\infty} B(K_m,\bB_m).
\end{equation}
\vskip .2cm
\begin{lem}\label{fucksian}
Let $B(K,\bB_0)$ be a $k$-beaded necklace subordinated
to a smooth $n$-knot $K$. Let $\Gamma_{\bB_0}$ be the group
generated by inversions through each $(n+1)$-sphere $\Sigma^{n+1}_{j}=\partial B^{n+2}_{j}$, where $B^{n+2}_{j}\in \bB_0$, $j=1,2,\ldots,k$.
Then the inverse limit space $\mathcal{K}$ is an $n$-knot embedded in $\mathbb{S}^{n+2}$.
\end{lem}

\noindent {\it Proof.} We will prove that $\mathcal{K}$ is homeomorphic to $\mathbb{S}^n$ by comparing our limit set $\mathcal{K}$ with a  well-known model (compare \cite{maskit} and \cite{DH} for $n=1$).\\

\noindent Let $B(K,\bB_0)=K\cup{\bB_0}$ be our beaded necklace, where  $\bB_0=\overset{k}{\underset{i=1}\bigcup}\,B^{n+2}_i$ is the finite union of disjoint, closed $(n+2)$-balls $B^{n+2}_i$ ($i\in\lbrace1,\dots,k\rbrace, k\geq3$). Now we will construct our beaded necklace model, consider the unit $n$-sphere $\mathbb{S}^{n}=\{x=(x_1,x_2,\ldots,x_{n+1},0)\in \mathbb{R}^{n+2} \,:\, \parallel x\parallel=1\}$ and we place $k$ disjoint closed Euclidean $(n+2)$-balls $O^{n+2}_i$ of radius $r>0$ centered on the points $c_i\in \mathbb{S}^{n}$. We also require that each ball be orthogonal to $\mathbb{S}^{n}$.  Thus we have ${\bO}_0=\bigcup_{i=1}^k O^{n+2}_i$ and the corresponding \emph{trivial $k$-beaded necklace} is given by $B(\mathbb{S}^n,{\bO}_0)=\mathbb{S}^n\cup{\bO}_{0}$.\\

\begin{figure}[ht] 
 \begin{center}
 \includegraphics[height=5cm]{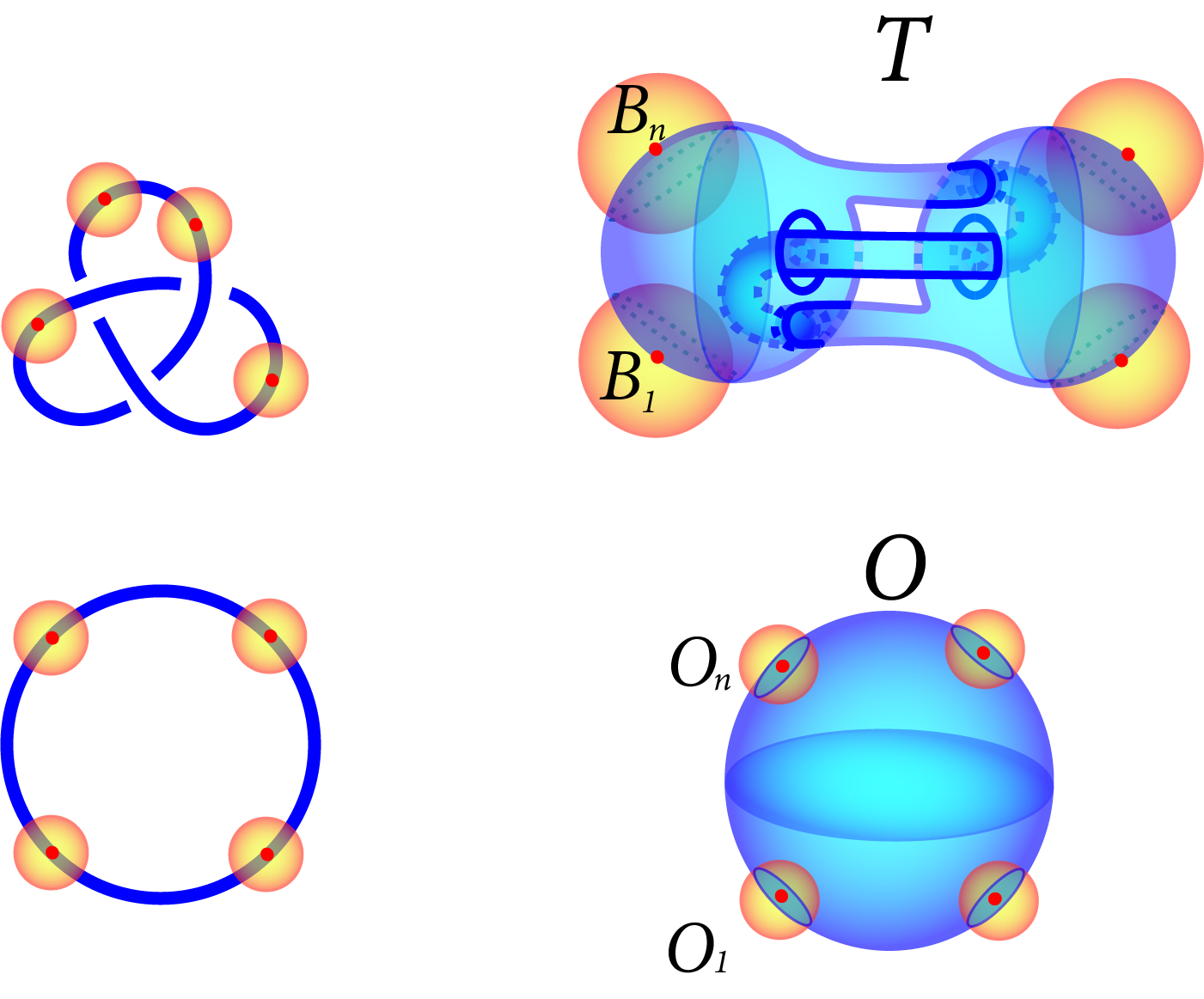}
\end{center}
\caption{\sl Schematic picture of the $k$-beaded necklaces $B(K,\bB_0)$ and $B(\mathbb{S}^n,{\mathcal{O}}_0)$.} 
\label{F15}
\end{figure}

\noindent Observe that, there exists a homeomorphism $h:B(K,\bB_0)\rightarrow B(\mathbb{S}^n,{\bO}_0)$  sending $B^{n+2}_j$ onto $O^{n+2}_j$, since we can renumber if necessary (see Figure \ref{F15}).\\

\noindent Now we will compare our two beaded necklaces at each stage of the inversion process. So in the first stage, we apply the inversions $I_j$ and $R_j$ on  $\partial B^{n+2}_j$ and $\partial O^{n+2}_j$ ($j=1,2,\ldots,k)$ respectively. Notice that 
$I_j(B(K,\bB_0)\smallsetminus{B^{n+2}_{j}})$ is homeomorphic to 
$R_j(B(\mathbb{S}^n,{\bO}_0)\smallsetminus{O^{n+2}_{j}})$, and since the intersection subset  $\partial B^{n+2}_j\cap K$ goes onto the corresponding intersection subset $\partial O^{n+2}_j\cap \mathbb{S}^n$, then there is a homeomorphism 
$h_{1}:B(K_1,\bB_1)\rightarrow B(\mathbb{S}^n,{\bO}_1)$ such that the following diagram commutes
\begin{equation*}
\xymatrix{
B(K_1,\bB_1)\ \ar@{^(->}[r]\ar[d]_{h_1} \hspace{.1cm}& \hspace{.1cm} B(K,\bB_0)\ar[d]^{h}\\
B(\mathbb{S}^n,{\bO}_1) \ar@{^(->}[r] \hspace{.1cm}& \hspace{.1cm}B(\mathbb{S}^n,{\bO}_0)
}
\end{equation*}
where the row maps are the corresponding inclusion maps. Notice that $B(\mathbb{S}^n,{\bO}_1)$ is subordinated to the unit $n$-sphere $\mathbb{S}^n$.\\

\noindent We continue in this way, so in the $m^{th}$-stage  we again apply the corresponding inversions $I_j$ and $R_j$, in such a way  that 
$I_j(B(K_{m-1},\bB_{m-1})\smallsetminus{B^{n+2}_{j}})$ is homeomorphic to $R_j(B(\mathbb{S}^n,{\bO}_{m-1})\smallsetminus{O^{n+2}_{j}})$, and again the intersection subset 
$\partial B^{n+2}_j\cap K_{m-1}$ is sent onto the corresponding intersection 
subset $\partial O^{n+2}_j\cap \mathbb{S}^n$. As a consequence,  we have a homeomorphism 
$h_{m}:B(K_m,\bB_m)\rightarrow B(\mathbb{S}^n,{\bO}_m)$  such that the following diagram commutes
\begin{equation*}
\xymatrix{
B(K_m,\bB_m)\ \ar@{^(->}[r]\ar[d]_{h_m} \hspace{.1cm}& \hspace{.1cm} B(K_{m-1},\bB_{m-1})\ar[d]^{h_{m-1}}\\
B(\mathbb{S}^n,{\bO}_m) \ar@{^(->}[r] \hspace{.1cm}& \hspace{.1cm}B(\mathbb{S}^n,{\bO}_{m-1})
}
\end{equation*}
where the row maps are the corresponding inclusion maps. Notice that the beaded necklace $B(\mathbb{S}^n,{\bO}_m)$ is subordinated to the unit $n$-sphere 
$\mathbb{S}^n$, for  each 
$m\in\mathbb{N}$.\\

\noindent Summarizing, we have the following commutative diagram
\begin{equation*}
\xymatrix{
 B(K,\bB_0)\ar[d]_{h} \hspace{.1cm}&\hspace{.1cm} \cdots\hspace{.1cm}\ar@{_(->}[l]\ar[d]\hspace{.1cm}&\hspace{.1cm} B(K_m,\bB_m)\ar@{_(->}[l]\ar[d]_{h_m} \hspace{.1cm}&\hspace{.1cm}\cdots\ar@{_(->}[l]\ar[d] \hspace{.1cm}&\hspace{.1cm}\mathcal{K}\ar@{_(->}[l]\ar[d]\\
B(\mathbb{S}^n,{\bO}_0)\hspace{.1cm}&\hspace{.1cm}\cdots \hspace{.1cm}\ar@{_(->}[l]\hspace{.1cm}&\hspace{.1cm} B(\mathbb{S}^n,{\bO}_m)\ar@{_(->}[l]\hspace{.1cm}&\hspace{.1cm}\cdots\ar@{_(->}[l]\hspace{.1cm}& \hspace{.1cm}\mathcal{O}\ar@{_(->}[l]
}
\end{equation*}
\noindent where $\mathcal{O}$ denote the inverse limit space given by  
$$
\mathcal{O}=\varprojlim_{m} B(\mathbb{S}^n,{\bO}_m)=\bigcap_{m=0}^{\infty} B(\mathbb{S}^n,{\bO}_0).
$$

\noindent We have, by the universal property of the inverse limit, a continuous map $H: \mathcal{K}\rightarrow \mathcal{O}$, and again by the same argument,  there is a continuous map $G: \mathcal{O}\rightarrow \mathcal{K}$.  Since each $h_m$ is a homeomorphism, it follows that $G$ is the inverse function of $H$; hence
$$
\mathcal{K}:= \varprojlim_{m} B(K_m,\bB_m)\cong \mathcal{O}.
$$
Observe that  $\mathcal{O}=\mathbb{S}^{n}$, since each beaded necklace 
$B(\mathbb{S}^n,\bO_m)$ is subordinated to  $\mathbb{S}^n$ (for more details see \cite{maskit} and \cite{DH}), therefore
$$
\mathcal{K}\cong \mathbb{S}^{n} 
$$\hfill $\square$\\

\noindent  We would like to remark that, by construction, the limit set $\Lambda(\Gamma_{\bB_0})$ of the Kleinian group $\Gamma_{\bB_0}$ is contained in $\mathcal{K}$.

\begin{lem}\label{knot}
Let $B(K,\bB_0)$ be a $k$-beaded necklace subordinated
to a smooth $n$-knot $K$. Let $\Gamma_{\bB_0}$ be the group
generated by inversions through each $(n+1)$-sphere $\Sigma^{n+1}_{j}=\partial B^{n+2}_{j}$, where $B^{n+2}_{j}\in \bB_0$, $j=1,2,\ldots,k$.
Then the Kleinian limit set $\Lambda(\Gamma_{\bB_0})$ of the group $\Gamma_{\bB_0}$  is homeomorphic to a Cantor
set embedded in $\mathbb{S}^{n+2}$.
\end{lem}

\noindent {\it Proof.} By definition of the Kleinian limit set (\cite{maskit}), we have that $\Lambda(\Gamma_{\bB_0})=\bigcap_m^{\infty} \bB _m$; \emph{i.e.}, the Kleinian limit set is the intersection of a nested sequence $\{\bB_m\}_{m\in\mathbb{N}}$, where each $\bB_m$ consists of the union of $l_m$ disjoint closed $(n+2)$-balls. 
Therefore, a point $x$ belongs to $\Lambda(\Gamma_{\bB_0})$ if and only if there exists a nested sequence of $(n+2)$-balls 
$\{B^{m}_{i_m}\}_{m\in\mathbb{N}}$, such that $x=\bigcap_m B^{m}_{i_m}$. Therefore, the result follows. $\square$

\begin{lem}\label{wildpoints}
Let $B(K,\bB_0)$ be a $k$-beaded necklace subordinated
to the non-trivial, smooth $n$-knot $K$ for $n>1$ \emph{with fundamental group {\bf bigger} than $\Z$}. Let $\Gamma_{\bB_0}$ be the group
generated by inversions through each $(n+1)$-sphere $\Sigma^{n+1}_{j}=\partial B^{n+2}_{j}$, where $B^{n+2}_{j}\subset \bB_0$, $j=1,2,\ldots,k$, whose Kleinian limit set is $\Lambda(\Gamma_{\bB_0})$. Then the corresponding inverse limit space $\mathcal{K}$ is wild at every point of $\Lambda(\Gamma_{\bB_0})$.
\end{lem} 

\noindent {\it Proof.} (Compare \cite{Hin} and \cite{DH}). Let $x\in
\Lambda(\Gamma_{\bB_0})$ be a limit point. By construction, $x\in\cK$. We want to show that $x$ is a wild point of $\cK$. Given an open 
neighborhood $U$ of $x$, by Equation \ref{beadsintersection}, there exists $m\geq1$ such
that $x$ is in the beaded necklace $B(K_m,\bB_m)$; therefore $x$ is contained in the interior one of the beads (a round closed ball of dimension $n+2$) of $B(K_m,\bB_m)$. Let's call this ball $B_m$.
Note that by construction, there are infinitely many copies of $K$ and its mirror image 
$\overline{K}$ in $B_m$.  Since by Equation  \ref{beadsintersection} one has that $\cK$ is the intersection of the bead necklaces $B(K_n,\bB_n)$ it follows that
\[
\bigcap_{n\geq{m}}B_m\cap{B(K_n,\bB_n)}=B_m\cap\cK
\]
\noindent Now, considering the set of complements one has:
\[
\bigcup_{n\geq{m}}B_m{\smallsetminus}{B(K_n,\bB_n)}=
B_m{\smallsetminus}\cap\cK
\]
\noindent Therefore, $\pi_1(B_m\smallsetminus\cK)$ is the direct limit
of $G_n\overset{def}=\pi_1(B_m{\smallsetminus}B(K_n,\bB_n))$:
\[
\underset{\iota_n}{\varinjlim}\,\left[\iota_n:G_n\to{G_{n+1}}, \,\, n\geq{m}\right],
\]
\noindent where $\iota_n$ is the homomorphism induced by the inclusion
 of $B_m{\smallsetminus}B(K_n,\bB_n))$ into $B_m{\smallsetminus}B(K_{n+1},\bB_{n+1})$.
 Lemma \ref{collars-retract} implies that the pair $(B_m, B_m\cap{B(K_{n},\bB_{n})})$
has the homotopy type of the tangle $(B_m, B_m\cap{K_n})$.
  This tangle is the tangle
associated with a knot of the form $K^{\# r_1(n)}\#(\overline{K})^{\# r_2(n)}\,,\,
 \, r_1(n), r_2(n)\in \N$;
where $\overline{K}$  is the mirror image of $K$.  
Therefore, using van Kampen's theorem 
(Proposition \ref{vankampen}, Corollary \ref{conn} and Remark \ref{generatorsgrow}), 
we obtain: 
\[
G_n\cong\;
\underbrace{G\;*_{\langle \mu\rangle}\;G\;*_{\langle \mu\rangle}\;\cdots\;*_{\langle \mu\rangle}\;G}_{\text{$(r_1(n)+r_2(n))$ copies}}
\]
%%%%%%%
%%%%%%%
where $G$ is the fundamental group of the canonical tangle associated to $K$ (Equation \ref{canonicaltangle}) and $\mu$ is a representative of the meridian (Definition \ref{ntangle}).\\

\noindent Let $G$ have the presentation
\[
G = \langle \mu, g_1, \dots, g_k \mid r_1, \dots, r_l\rangle,
\]
where $\mu$ is the meridian loop around the knot. 
Since we are assuming that the fundamental group of $\bS^{n+2}{\smallsetminus}K$  is bigger than $\Z$ we must have that $k\geq1$ and the minimal cardinality of a set
of generators of $G_n$ must be at least $r_1(n)+r_2(n)$, as $n\mapsto\infty$ 
$r_1(n)+r_2(n)\mapsto\infty$.\\

%%%%%%%
%%%%%%%
\noindent Summarizing:
\[
\pi_1\left(B_m\setminus \cK)\right) \cong \left( \ast_{i=1}^\infty G_i \right) \Big/ \left\langle \mu_i = \mu_{i+1} \ \forall i \right\rangle,
\]
where each $G_i\cong{G}$ is bigger than $\Z$, 
$\mu_i$ is a representative of the meridian in $G_i$, and the identification of meridians across copies does not reduce the group to a finitely generated one.\\

\noindent Since the fundamental group of $K$ is bigger than $\Z$, this implies that
$\pi_{1}(B_m\smallsetminus\mathcal{K})$ is not isomorphic 
to a finitely generated group, \ie it
is infinitely generated. Since we can choose $B_m$ arbitrarily small, the knot must be wild, and the result follows. $\square$

\begin{rem} 
There are $n$-knots $K^n\subset{S^{n+2}}$, for $n\geq3$ such that 
$\pi_1(S^{n+2}\smallsetminus{K^n})=\Z$, for instance, the examples given by Milnor´s Fibration Theorem given below in example \ref{ex3} in Examples \ref{example-fibered}.
However, for any dimension $n\geq3$ there are examples of knots with fundamental groups bigger than $\Z$. For instance, spin knots. In fact, a theorem by Kervaire \cite{ker1, ker2} states that if  $n \geq 5$, a finitely presented group $G$ is the group of an $n$-knot 
if and only if it satisfies the following conditions:
\begin{enumerate}
    \item $H_1(G,\Z) \cong \mathbb{Z}$,
    \item $H_2(G,\Z) = 0$,
    \item $G$ is normally generated by a single element.
    \item $G$ is the fundamental group of a smooth homology $(n+2)$-sphere. 
        That is, there exists a smooth, closed, oriented $(n+2)$-manifold $M$ 
        with $H_*(M) \cong H_*(S^{n+2})$ and $\pi_1(M) \cong G$.
   \end{enumerate}
   \noindent For any $n\geq5$, there exist groups satisfying this condition.
\end{rem}

\begin{definition}
The wild $n$-knot $\mathcal{K}$ is called a \emph{wild $n$-knot invariant by a Schottky group}.
\end{definition}

\begin{rem}\label{addr}
  By the previous discussion, we have that given $a\in \Lambda(\Gamma_{\bB_0})$, there exists a sequence
  $\{B^{l}_{i_l}\}_{l\in\mathbb{N}}$
  such that $a\in B^{l}_{i_l}$ and $B^{l}_{i_l}\subset B(K_l,\bB_l)$. Since 
  $\{B^{l}_{i_l}\}_{l\in \mathbb{N}}$ is a sequence of nested sets and
  $\diam(B^{l}_{i_l})\to 0$, by
  Cantor's IntersectionTheorem, it follows that $\{a\} =\bigcap B^{l}_{i_l}$. Therefore $\{i_0,\,i_1,\dots,i_l,\dots\}$
  is an address
  for each $a\in \Lambda(\Gamma_{\bB_0})$ (see Figure \ref{Schem} and compare \cite{GCA}).  
\end{rem}

\begin{rem}
We recall that an $n$-knot $K$ is homogeneous, if  given two points
$p,\,\,q\in{K}$ there exists a homeomorphism $f:\mathbb{S}^{n+2}\rightarrow\mathbb{S}^{n+2}$
such that $f(K)=K$ and $f(p)=q$. Observe that by our construction, any wild $n$-knot invariant by a Schottky group is not homogeneous, since it contains both wild and tame points. However, there exist homogeneous wild $1$-knots, for instance, any dynamically defined wild 1-knot is homogeneous (for more details see \cite{GA} and \cite{BHV}).
\end{rem}

\section{Fibration of $\mathbb{S}^{n+2}\smallsetminus\mathcal{K}$ over
 $\mathbb{S}^{1}$}
 
We recall that an $n$-knot or $n$-link $L$ in $\mathbb{S}^{n+2}$ is {\it fibered} if there exists a 
locally trivial fibration $f:(\mathbb{S}^{n+2}\smallsetminus{L})\rightarrow \mathbb{S}^{1}$. We 
require that $f$ be well-behaved near $L$, that is, each component $L_{i}$ has a neighborhood framed as $\mathbb{D}^{2}\times\mathbb{S}^{n}$, with 
$L_{i}\cong \{0\}\times\mathbb{S}^{n}$, in such a way that the
restriction of $f$ to $(\mathbb{D}^{2}{\smallsetminus}\{0\})\times\mathbb{S}^{n}$ is the map 
into $\mathbb{S}^{1}$ given by $(x,y)\rightarrow \frac{x}{|x|}$. It follows that each
$f^{-1}(x)\cup L$, $x\in\mathbb{S}^{1}$, is an $(n+1)$-manifold
 with boundary $L$; in fact a Seifert ``surface'' for $L$ (see \cite{rolfsen}, page 323).

\begin{main1}
Let $B(K,\bB_0)$ be a $k$-beaded necklace subordinated
to the non-trivial, smooth fibered $n$-knot $K$. Let $\Gamma_{\bB_0}$ be the group
generated by inversions through each $(n+1)$-sphere $\Sigma^{n+1}_{j}=\partial B^{n+2}_{j}$, where $B^{n+2}_{j}\subset \bB_0$, $j=1,2,\ldots,k$, and consider the 
inverse limit space  $\mathcal{K}$.  Then:
\begin{enumerate}
\item There exists a locally trivial fibration $\psi
:(\mathbb{S}^{n+2}\smallsetminus\mathcal{K})\rightarrow\mathbb{S}^{1}$, where the
 fiber $\Sigma_{\theta}=\psi^{-1}(\theta)$ is an orientable $(n+1)$-manifold with one end,
and if the fiber of $K$ has non-trivial homology in dimension $r$,
then $H_r(\psi^{-1}(\theta),\Z)$ is infinitely generated.
\item  $\overline{\Sigma_{\theta}}\smallsetminus\Sigma_{\theta}=\mathcal{K}$, where
$\overline{\Sigma_{\theta}}$ is the closure of $\Sigma_{\theta}$ in $\mathbb{S}^{n+2}$.
\end{enumerate}
\end{main1}

\begin{figure}[ht] 
 \begin{center}
  \includegraphics[height=3cm]{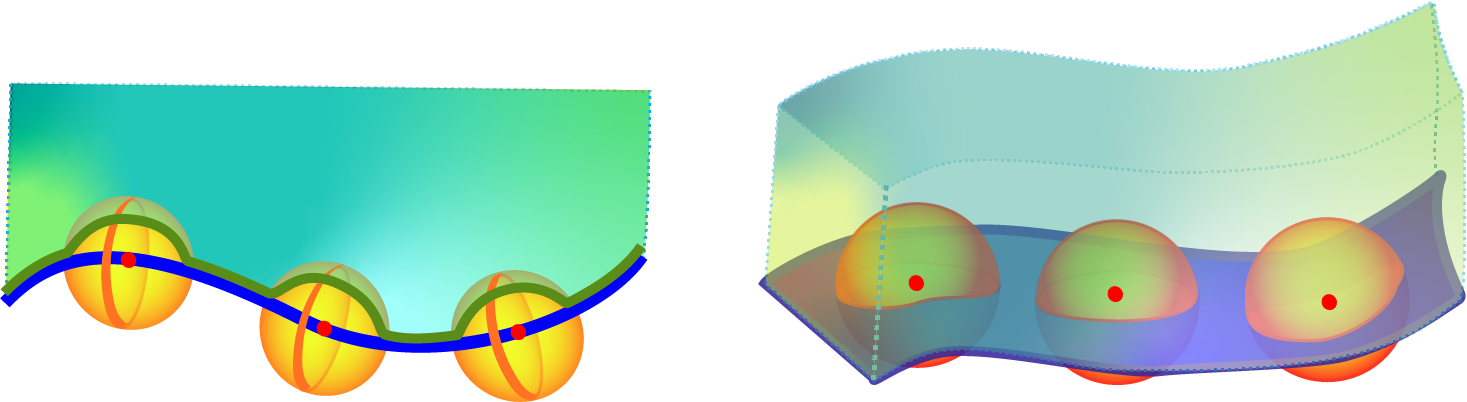}
\end{center}
\caption{\sl The fiber $S$, up to isotopy.} 
\label{F4}
\end{figure}

\noindent {\it Proof.} The proof is similar to the one given in \cite{DH} and \cite{Hin}.
We will prove that $\mathbb{S}^{n+2}\smallsetminus\mathcal{K}$ fibers over the circle. Since $K$ fibers over the circle, consider the corresponding fibration
$\widetilde{P}:(\mathbb{S}^{n+2}\smallsetminus{K})\rightarrow \mathbb{S}^{1}$ with fiber the oriented $(n+1)$-manifold $S$. Then the restriction map 
$\widetilde{P}\mid_{\mathbb{S}^{n+2}\smallsetminus{B(K,\bB_0)}}\equiv
P$ is also a fibration and, after modifying $\widetilde{P}$ by isotopy if necessary, we can consider  that
the fiber $S$  intersects the boundary of each bead $B^{n+2}_{i}\in \bB_0$ in an $n$-disk
$a_{i}$, whose boundary belong to $\mathcal{K}$ (see Figure \ref{F4}). The 
fiber $\widetilde{P}^{-1}(\theta)=P^{-1}(\theta)$ is the $(n+1)$-manifold 
$S^{*}=\bar{S}\smallsetminus{K}$ for any
$\theta\in\mathbb{S}^{1}$, where $\bar{S}$ denotes the closure of the
$(n+1)$-manifold $S$ in $\mathbb{S}^{n+2}$. Thus $S^{*}$ is oriented and its boundary  intersects
each $(n+1)$-sphere $\Sigma_{j}$ in an $n$-disk $a_{j}$ whose boundary $\partial a_{j}$ is contained into $K$.\\

\begin{figure}[ht] 
 \begin{center}
 \includegraphics[height=3cm]{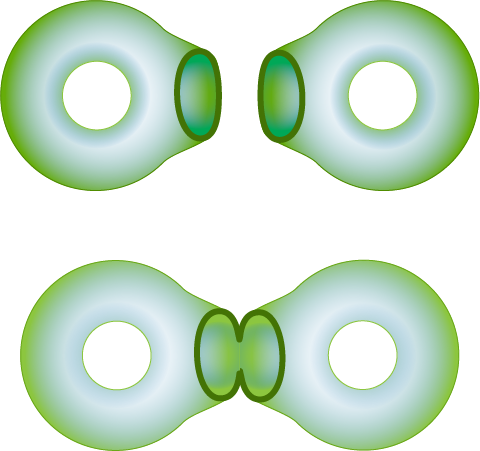}
\end{center}
\caption{\sl Sum of two $(n+1)$-manifolds $S^{*}$ and $S^{*1}_{j_0}$ along the $n$-disk $a_{j_0}$.} 
\label{F6}
\end{figure}

\noindent Observe that the inversion map  $I_{j_0}$ sends both a copy of 
$B(K,\bB_0)\smallsetminus{B^{n+2}_{j_0}}$
(namely $\kappa^1_{j_0}=I_{j_0}(B(K,\bB_0)\smallsetminus{B^{n+2}_{j_0}})$ and a copy of $S^{*}$ (called $S^{*1}_{j_0}$) into
the bead $B^{n+2}_{j_0}$, for $j_0=1,2,\ldots,k$. Since both $\kappa^1_{j_0}$ and $S^{*1}_{j_0}$ have opposite orientations.
Then the beaded necklaces $\kappa^1_{j_0}$ and $B(K,\bB_0)$ are joined by the set  $\partial B^{n+2}_{j_0}\cap K$, and  $S^{*}$ and $S^{*1}_{j_0}$ are joined by
the $n$-disk $a_{j_0}$ (see Figure \ref{F6}) which, in both manifolds, has the same orientation. \\

\noindent So, we have a new 
beaded necklace $B(K_{1_{j_0}},\bB_{1_{j_{0}}})=(B(K,\bB_0)\smallsetminus{B^{n+2}_{j_0}})\cup \kappa^1_{j_0}$ isotopic to the 
connected sum $B(K,\bB_0)\#\kappa^1_{j_0}$, whose complement also fibers over the circle with 
fiber the sum of $S^{*}$ with $S^{*1}_{j_0}$ along the $n$-disk $a_{j_0}$, namely
$S^{*}\#_{a_{j_0}}S^{*1}_{j_0}$. In fact, the map 
$P_{1_{j_0}}:(\mathbb{S}^{n+2}\smallsetminus{B(K_{1_{j_0}}},\bB_{1_{j_{0}}}))\rightarrow \mathbb{S}^1$ given by
$$
P_{1_{j_0}}(x)=\left\{ \begin{array}{ll}
            P(x) & \mbox{if $x \in\mathbb{S}^{n+2}\smallsetminus{B(K,\bB_{0}})$}\\
P I_{j_0}(x) &\mbox{if $ x\in B^{n+2}_{j_0} \smallsetminus{I_{j_0}(B(K,\bB_0)}\smallsetminus{B^{n+2}_{j_0}})$, where $j_0\in \{1,\ldots ,k\}.$} \end{array} \right.  
$$ 
is a well-defined locally trivial fibration  map, since both $P(x)$ and $P I_{j_0}$ are continuous and coincide at the intersection. The  fiber  is $S\#_{a_{j_0}}S^1_{j_0}$. Now, if the fiber $S^{*}$ has homology $H_r(S^{*},\Z)$ in dimension $r>0$, then using the Mayer-Vietoris sequence and the fact that $H_r(a_{j_0},\Z)$ is trivial for $r>0$, we have that 
$H_r(S\#_{a_{j_0}}S^1_{j_0},\Z)\cong H_r(S^{*},\Z)\oplus H_r(S^{*},\Z)$. \\

\begin{figure}[ht] 
 \begin{center}
 \includegraphics[height=3.5cm]{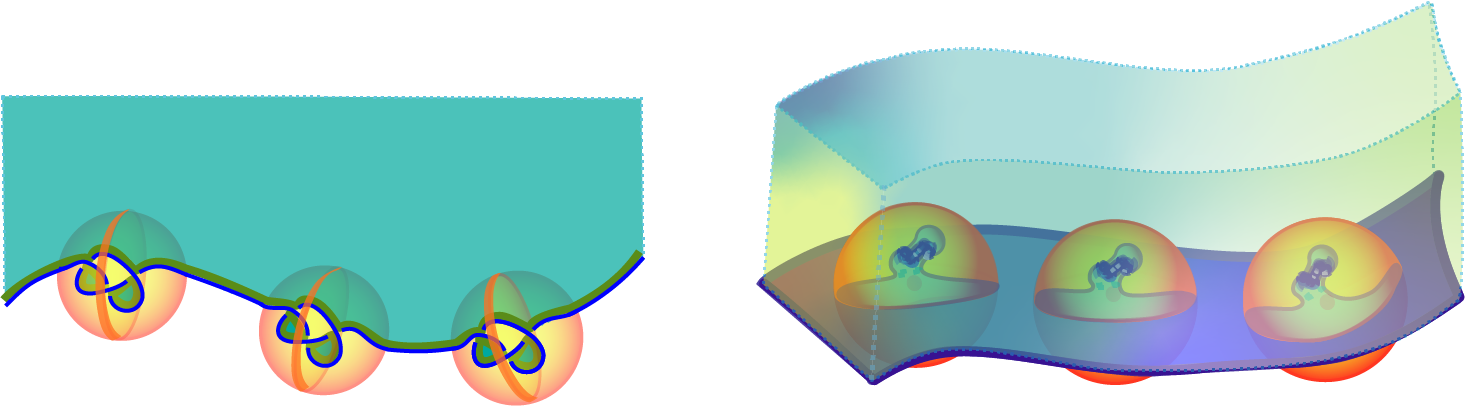}
\end{center}
\caption{\sl A schematic picture of the fiber $S^{*1}$.} 
\label{F7}
\end{figure}

\noindent At the end of the first stage, we get  a new beaded necklace $B(K_1,\bB_1)$ subordinated to the new knot $K_{1}$
(see Section \ref{dyn}) such that, by the previous discussion, its complement fibers 
over the circle via the locally trivial fibration map 
$P_1:\mathbb{S}^{n+2}\smallsetminus{B(K_1,\bB_1)}\rightarrow \mathbb{S}^1$ given by  
$$
P_{1}(x)=\left\{ \begin{array}{ll}
            P(x) & \mbox{if $x \in\mathbb{S}^{n+2}\smallsetminus{B(K,\bB_0})$}\\
P I_{j_0}(x) &\mbox{if $ x\in B^{n+2}_{j_0}\smallsetminus{I_{j_0}(B(K,\bB_0)}\smallsetminus{B^{n+2}_{j_0}})$, where $j_0\in \{1,\ldots ,k\}.$} \end{array} \right.  
$$ 
The fiber is the $(n+1)$-manifold $S^{*1}$ which is, in turn, 
homeomorphic to the sum of $k+1$ copies of $S^{*}$ along the respective $n$-disks and its homology groups are given by $H_r(S^{*1},\Z)\cong \bigoplus_{k+1} H_r(S^{*},\Z)$ for $r>0$ (see Figure \ref{F7}).\\

\noindent Then  the following diagram commutes
\begin{equation*}
\xymatrix{
\mathbb{S}^{n+2}\smallsetminus{B(K,\bB_0)} \ar@{^(->}[r]\ar[d]_{P} \hspace{.1cm}& \hspace{.1cm}\mathbb{S}^{n+2}\smallsetminus{B(K_1,\bB_1)}\ar[d]_{P_1}\\
\mathbb{S}^1\ar@{^(->}[r] \hspace{.1cm}& \hspace{.1cm}\mathbb{S}^1
}
\end{equation*}
where the bottom row map is the identity. \\

\noindent Continuing with the inverting process, at the $m^{th}$ stage we obtain 
$\bB_{m}$ which is the union of  $l_m=k(k-1)^{m}$ beads $B_j^{m}$, $j\in\{1,\dots, l_m\}$, subordinated to a tame 
$n$-knot $K_{m}$.  Let  $B(K_m,\bB_m)=K_m\cup \bB_{m}$ be the corresponding new beaded necklace, so by construction $\bB_{m}\subset \bB_{m-1}$ and 
$B(K_m,\bB_m)\subset B(K_{m-1},\bB_{m-1})$. Notice that the locally trivial fibration map $P_{m-1}:(\mathbb{S}^{n+2}{\smallsetminus}B(K_{m-1},\bB_{m-1}))\rightarrow \mathbb{S}^1$ with fiber the Seifert surface $S^{*m-1}$, can be extended to a continuous map 
$P_m:(\mathbb{S}^{n+2}{\smallsetminus}B(K_m,\bB_m))\rightarrow \mathbb{S}^1$, as follows
$$
P_{m}(x)=\left\{ \begin{array}{ll}
            P_{m-1}(x) & \mbox{if $x\in \mathbb{S}^{n+2}{\smallsetminus}B(K_{m-1},\bB_{m-1})$}\\
            P I_{j_0} \cdots I_{j_{m-1}}I_{j_m}(x) &\mbox{if $ x\in B^{n+2}_{j_0,\ldots, j_{m}}{\smallsetminus}$}\\
                                                                                                                &{\smallsetminus}I_{j_m}(B(K_{m-1},\bB_{m-1}){\smallsetminus}B^{n+2}_{j_0,\ldots, j_{m}}) \end{array} \right.  
$$ 
where $B^{n+2}_{j_0,\ldots, j_{m}}:=I_{j_m}I_{j_{m-1}}\ldots I_{j_1}(B^{n+2}_{j_0})$ for $B^{n+2}_{j_0}\subset \bB_0$, and $j_s\in\{1,2,\ldots,k\}$.  So, $P_m$ is a locally-trivial fibration map, since  $P_{m-1}(x)$ and $P I_{j_0} \cdots I_{j_{m-1}}I_{j_m}$ are continuous and coincide in the intersection. The fiber is the $(n+1)$-manifold $S^{*m}$ which is, in turn, homeomorphic to the sum of $l_m+1$ copies of $S^{*}$ along the respective 
$n$-disks  (see Figure \ref{F7}). As above,  we have that the homology groups of $S^{*m}$ are $H_r(S^{*m},\Z)\cong \bigoplus_{l_{m+1}} H_r(S^{*},\Z)$ for $r>0$.
Notice that in each step, the boundaries of the $n$-disks are removed, since they are contained in $\mathcal{K}$,
and the diameter of  $a_{j}$ tends to zero.  Then the following diagram commutes:\\

\begin{equation*}
\xymatrix{
\mathbb{S}^{n+2}{\smallsetminus}B(K_{m-1},\bB_{m-1}) \ar@{^(->}[r]\ar[d]_{P_{m-1}} \hspace{.1cm}& \hspace{.1cm}\ar[d]_{P_m}\mathbb{S}^{n+2}{\smallsetminus}B(K_m,\bB_m)\\
\mathbb{S}^1\ar@{^(->}[r] \hspace{.1cm}& \hspace{.1cm}\mathbb{S}^1
}
\end{equation*}
where the bottom row map is the identity. \\

\noindent Summarizing, we have the following commutative diagram:\\

\begin{equation*}
\xymatrix{
\mathbb{S}^{n+2}{\smallsetminus}B(K,\bB_0)\ar@{^(->}[r]\ar[d]_{P} \hspace{.1cm}&\hspace{.1cm}\mathbb{S}^{n+2}{\smallsetminus}B(K_1,\bB_1)\ar@{^(->}[r]\ar[d]_{P_1} \hspace{.1cm}&\hspace{.1cm}\cdots\ar@{^(->}[r]\ar[d] \hspace{.1cm}&\hspace{.1cm}\mathbb{S}^{n+2}{\smallsetminus}\mathcal{K}\ar[d]\\
\mathbb{S}^1\ar@{^(->}[r]\hspace{.1cm}&\hspace{.1cm}\mathbb{S}^1\ar@{^(->}[r]\hspace{.1cm}&\hspace{.1cm}\cdots\ar@{^(->}[r]\hspace{.1cm}&\hspace{.1cm}\mathbb{S}^1
}
\end{equation*}

\noindent By the universal property of the direct limit, there exists a continuous function $\psi:(\mathbb{S}^{n+2}{\smallsetminus}\mathcal{K})\rightarrow
\mathbb{S}^1$. Since each $P_m$ is onto, it follows that $\psi$ is also onto. \\

\noindent Given $\theta\in \mathbb{S}^1$. For each integer $m$, 
$P_m^{-1}(\theta)$ is the $(n+1)$-manifold $S^{*m}$, and since $P_m$ is a locally trivial fibration map, there exists an open neighborhood $W_m$ of $\theta$ in $\mathbb{S}^1$ such that $P^{-1}_m(W_m)$ is homeomorphic to $W_m\times S^{*m}$. Even more, by construction, we have that  $W_m:=W_0$ for each $m$, hence 
 $\psi^{-1}(W_0)$ is homeomorphic to $W_0\times \Sigma_{\theta}$, where $\psi^{-1}(\theta):=\Sigma_{\theta}$ is the fiber which is, the direct limit of 
$\{S^{*m},\hspace{.2cm}m=0,1,\ldots\,:\hspace{.2cm}i_{k}:S^{*k}\rightarrow S^{*k+1}\}$, where $i_{k}$ is the inclusion map.
It is the sum along $n$-disks of an infinite
number of copies of $S^{*}$, hence it is homeomorphic to an orientable $(n+1)$-manifold. This implies that the homology groups of $\Sigma_{\theta}$ are the direct limit of  $\{H_r(S^{*m},\Z),\hspace{.2cm} m=0,1,\ldots\,:\hspace{.2cm}i_{*m}:H_r(S^{*m},\Z)\rightarrow H_r(S^{*m+1},\Z)\}$, where $i_{*k}$ is the corresponding inclusion map, for $r>0$. As a consequence, if $S$ has non-trivial homology in dimension $r$
then $H_r(\psi^{-1}(\theta),\Z)$ is infinitely generated. Now, we will describe its set of ends.
\noindent Consider the Fuchsian model (see \cite{maskit}). So, the beaded necklace $B(\mathbb{S}^n,{\bO}_0)=\mathbb{S}^n\cup {\bO}_{0}$ consists of the unit $n$-sphere and  ${\bO}_0=\bigcup_{i=1}^k O^{n+2}_i$ where  $O^{n+2}_i$ is a closed Euclidean $(n+2)$-ball of radius $r>0$ and centered on points $c_i\in \mathbb{S}^{n}$. We also require that each ball be orthogonal to $\mathbb{S}^{n}$. Then its limit set
is the $n$-sphere and its complement fibers over $\mathbb{S}^{1}$
with fiber the $(n+1)$-disk  (for more details see the proof of Lemma \ref{fucksian}). \\

\noindent Returning to our case, if we intersect this disk with any
compact set, its complement consists of just one connected component. Hence, it has only
one end. Therefore, our manifold has one end.\\

\noindent The first part of the theorem has been proved. For the second part,
observe that the closure of the fiber in $\mathbb{S}^{n+2}$ is the fiber
union of its boundary. Therefore 
$\overline{\Sigma_{\theta}}{\smallsetminus}\Sigma_{\theta}=\mathcal{K}$. 
\begin{flushright}
$\square$
\end{flushright} 
%\vskip .2cm
\begin{rem}
Let $K$ be a non-trivial fibered smooth $n$-knot whose fiber is the oriented $(n+1)$-manifold $S$. So  $\mathbb{S}^{n+2}{\smallsetminus}K$ admits an open book decomposition, where each page is $S$ and its binding is $K$. We can describe the monodromy of $\mathbb{S}^{n+2}{\smallsetminus}K$  via the first return Poincar\'e map $\Phi$ that
 is the flow that cuts transversally each page of its open book decomposition (see \cite{verjovsky}, chapter 5). \\
 
\noindent As in wild 1-knots of dynamically defined types' case (see \cite{DH}), consider a beaded necklace $B(K,\bB_0)$ subordinated to $K$, then throughout the inverting process, $K$ and $S$ are copied into each ball by the corresponding inversion map (preserving or reversing orientation). So the flow $\Phi$ is also copied, and its direction changes according to the number of inversions, hence the 
Poincar\'e map can be extended at each stage, providing us in the limit a
homeomorphism $\psi:\Sigma_{\theta}\rightarrow \Sigma_{\theta}$ that identifies 
$\Sigma_{\theta}\times\{0\}$ with $\Sigma_{\theta}\times\{1\}$ and which induces
the monodromy of the corresponding wild $n$-knot.\\

\noindent Therefore the {\it monodromy} is an {\it invariant} for  wild $n$-knots invariant by a Schottky group.
\end{rem}

\begin{exs}\label{example-fibered}
In the following examples, we can explicitly describe the fundamental group of the complement as the semidirect product of the integer $\mathbb{Z}$ with the infinite free product of the fundamental group of the fiber of the original fibered knot with itself.

\begin{enumerate}
\item  \label{ex1} The unknot $\mathbb{S}^n\subset\mathbb{S}^{n+2}$ is fibered by the projection map $(\mathbb{S}^n\star\mathbb{S}^1){\smallsetminus}\mathbb{S}^n\rightarrow\mathbb{S}^1$, where $\star$ represents the join of spaces. In this case, the fibers are $(n+1)$-disks.

\item  \label{ex2} (Cappell-Shaneson \cite{CS}). Let $A$ be a Cappell-Shaneson $3\times3$-matrix, \emph{i.e,} $A\in\text{SL}(3,\Z)$ and $\text{det}(A)=1$. Then $A$ induces an automorphism $A:\bT^3\to\bT^3$, of the 3-torus (as an abelian Lie group with identity $e$). Let $N$ be the mapping torus of $A$. Then $N$ is a 4-manifold with fibre $T^3$ and monodromy $A$. Then
$N=\bT^3\times[0,1]/\sim$ where $(x,0)\sim(A(x),1)$ in $\bT^3\times[0,1]$. 
Since $A(e)=e$, under this identification, $e\times[0,1]$ becomes a circle $C\subset{N}$
with a trivial normal bundle in $N$. Using surgery, one replaces the interior of a tubular neighborhood $B^3\times{C}=B^3\times\bS^1$ by the interior of $\bS^2\times{B^2}$. 
One can show \cite{CS} that the resulting 4-manifold $\hat{N}$ is homeomorphic to the 4-sphere $\bS^4$ (for infinitely many choices of $A$, $N$ is the 
standard 4-sphere, \ie it is diffeomorphic to the 4-sphere with its usual differentiable structure \cite{Ak}). Then $K=\bS^2\times{0}$ becomes a 2-dimensional 
fibered knot in $\hat{N}$, with fibre $\bT^3{\smallsetminus}\lbrace\text{e}\rbrace$.

\item  \label{ex3} Milnor's fibration theorem provides beautiful examples of high-dimensional fibered knots \cite{Mi}. For instance, Brieskorn proved in \cite{B} that if $f:\C^n\to\C$, $f(z_1,\dots,z_n)=z_1^3+z_2^2+\dots+z_n^2$, then $f^{-1}(0)\cap\bS^{2n-1}$ is an exotic sphere (\cite{B},\cite{KM}, \cite{Mi}) if $n\geq5$. Furthermore, the complement fibers over $\bS^1$ and the fiber has the homotopy type of a wedge of two
$(n-1)$-spheres. Then the fibre $F$ of the associated wild knot has that $H_{n-1}(F,\Z)$
is infinitely generated. However, for $n\geq5$ the complement of the knot has a fundamental group isomorphic to $\Z$\, so that the associated knot also has a fundamental group isomorphic to $\Z$.

\item The right-handed trefoil 1-knot and the figure-eight 1-knot are fibered knots with fiber
the punctured torus.
\end{enumerate}
\end{exs}

\section{Equivalent wild $n$-knots invariant by a Schottky group}

In this section, we will prove Theorem 2.

\begin{definition}\label{NE}
Let $B(K,\bB_0)$ and $B(L,\bC_0)$ be two $k$-beaded necklaces subordinated to the smooth $n$-knots $K$ and $L$ respectively, such that $\bB_0=\bigcup_{j=1}^k B_j$ and $\bC_0=\bigcup_{i=1}^k C_i$ where $B_j$, $C_i$ are closed $(n+2)$- balls $(i,\,j=1,2,\dots,k)$. We say that $B(L,\bC_0)$ is 
\emph{equivalent} to $B(K,\bB_0)$, if  there exists a homeomorphism $\varphi:\mathbb{S}^{n+2}\rightarrow \mathbb{S}^{n+2}$ such that 
$\varphi (B(K,\bB_0))=B(L,\bC_0)$ and $\varphi(B_i)=C_i$.
\end{definition}

\begin{rem}
Let $B(K,\bB_0)$ and $B(L,\bC_0)$ be two equivalent $k$-beaded necklaces subordinated to the smooth $n$-knots $K$ and $L$, respectively. Then $K$ and $L$ are equivalent $n$-knots.
\end{rem}

\begin{main2}\label{main2}
Let $B(K,\bB_0)$ and $B(L,\bC_0)$ be two equivalent $k$-beaded necklaces subordinated to the smooth $n$-knots $K$ and $L$, respectively. Then the corresponding inverse limit spaces $\mathcal{K}$ and $\mathcal{L}$ are equivalent wild $n$-knots.
\end{main2}

\noindent {\it Proof.}  Let $B(K,\bB_0)$ and $B(L,\bC_0)$ be two equivalent $k$-beaded necklaces via the homeomorphism $H:\mathbb{S}^{n+2}\rightarrow \mathbb{S}^{n+2}$.  Consider  $\bB_0=\bigcup_{i=1}^k B_i$ and $\bC_0=\bigcup_{i=1}^k C_i$ and let $\Gamma_{\bB_0}$ and $\Gamma_{\bC_0}$
 be the corresponding groups generated by inversions $I_{j}$, $J_j$ through $\Sigma_{j}=\partial B_j$ and $\Upsilon_j=\partial C_j$  ($j=1,\ldots,k$), respectively.\\

\noindent On the first stage of the inverting process, we obtain the $l_1$-beaded necklaces $B(K_1,\bB_1)$ and $B(L_1,\bC_1)$ from $B(K,\bB_0)$ and $B(L,\bC_0)$, respectively, where 
$l_1=k(k-1)$ (see Section \ref{dyn}). 
Notice that $B(K_1,\bB_1)$ and $B(L_1,\bC_1)$ are also equivalent beaded necklaces. Consider the function 
$H_1:\mathbb{S}^{n+2}\rightarrow \mathbb{S}^{n+2}$ given by 
$$
H_1(x)=\left\{ \begin{array}{ll}
            H(x) & \mbox{if $x \in\mathbb{S}^{n+2}{\smallsetminus}\bB_0$}\\
J_{j_1}H I_{j_1}(x) &\mbox{if $ x\in B^{n+2}_{j_1}{\smallsetminus}I_{j_1}(\bB_0 {\smallsetminus}B^{n+2}_{j_1})$, where $j_1\in \{1,\ldots ,k\}.$} \end{array} \right.  
$$ 
Then, it is a well-defined continuous function and since $I_{j}$, $J_j$ are  homeomorphisms ($j=1,\ldots,k$), we have $H_1$ is also a homeomorphism that satisfies  $H_1(B(K_1,\bB_1))=B(L_1,\bC_1)$. Hence,  the following diagram is commutative
\begin{equation*}
\xymatrix{
B(K_1,\bB_1)\ar[r]\ar[d]_{H_1} & B(K,\bB_0)\ar[d]^{H}\\
B(L_1,\bC_1)\ar[r] & B(L,\bC_0)
}
\end{equation*}
where the row maps are the corresponding inclusion maps. \\

\noindent Continuing with the inverting process, we have in the $m$-stage that the corresponding $l_m$-beaded necklaces $B(K_m,\bB_m)$ and 
$B(L_m,\bC_m)$  ($l_m=k(k-1)^m$) are equivalent via the homeomorphism $H_m:\mathbb{S}^{n+2}\rightarrow \mathbb{S}^{n+2}$ given by 
$$
H_m(x)=\left\{ \begin{array}{ll}
            H(x) & \mbox{if $x \in\mathbb{S}^{n+2}{\smallsetminus}\bB_m$}\\
J_{j_m}J_{j_{m-1}}\cdots J_{j_1}H I_{j_1}(x) \cdots I_{j_{m-1}}I_{j_m} &\mbox{if $ x\in B^{n+2}_{j_1,\ldots, j_{m}}$}\\
                                                                                                                &{\smallsetminus}I_{j_m}(\bB_{m-1}{\smallsetminus}B^{n+2}_{j_1,\ldots, j_{m}}) \end{array} \right.  
$$ 
where  $j_1,\ldots, j_m\in \{1,\ldots ,k\}$ (see Remark \ref{addr} and proof of Theorem 1). So $H_m$  is a well-defined continuous function and again since $I_{j}$, $J_j$ are  homeomorphisms ($j=1,\ldots,n$), it follows that $H_m$ is also a homeomorphism that satisfies  $H_m(B(K_m,\bB_m))=B(L_m,\bC_m)$. Hence, the  following diagram is commutative
\begin{equation*}
\xymatrix{
B(K_m,\bB_m)\ar[r]\ar[d]_{H_m} & B(K_{m-1},\bB_{m-1})\ar[d]^{H_{m-1}}\\
B(L_m,\bC_m)\ar[r] & B(L_{m-1},\bC_{m-1})
}
\end{equation*}
where the row maps again are the corresponding inclusion maps. \\

\noindent Summarizing, we have the following commutative diagram

$$
\begin{CD}
B(K,\bB_0)@<<<B(K_1,\bB_1)@<<<\cdots B(K_m,\bB_m)@<<<\cdots\mathcal{K}\\
@V H VV@V H_1 VV@V H_m VV@V VV\\
B(L,\bC_0)@<<<B(L_1,\bC_1)@<<<\cdots B(L_m,\bC_m)
@<<<\cdots \mathcal{L}\\
\end{CD}
$$
\vskip .2cm
\noindent By the universal property of the inverse limit, there exists a continuous function $\tilde{F}:\mathcal{K}\rightarrow \mathcal{L}$ and using the same argument, 
there exists a continuous function $\tilde{G}:\mathcal{L}\rightarrow \mathcal{K}$ and by the following commutative diagram, we have that $\tilde{F}\circ\tilde{G}=Id$ and 
$\tilde{G}\circ\tilde{F}=Id$.

$$
\begin{CD}
B(K,\bB_0)@<<<B(K_1,\bB_1)@<<<\cdots B(K_m,\bB_m)@<<<\cdots\mathcal{K}\\
@V H VV@V H_1 VV@V H_m VV@V VV\\
B(L,\bC_0)@<<<B(L_1,\bC_1)@<<<\cdots B(L_m,\bC_m)
@<<<\cdots \mathcal{L}\\
@V H^{-1} VV@V H^{-1}_1 VV@V H^{-1}_m VV@V VV\\
B(K,\bB_0)@<<<B(K_1,\bB_1)@<<<\cdots B(K_m,\bB_m)@<<<\cdots\mathcal{K}\
\end{CD}
$$
\vskip .2cm
\noindent Therefore $\mathcal{K}$ and $\mathcal{L}$ are equivalent wild $n$ knots invariant by a Schottky group. $\square$

\section{Cyclic branched covers of wild knots of dynamical type}

Consider a $k$-beaded necklace $B(K,\bB_0)$ subordinated
to the non-trivial, smooth $n$-knot $K$ for $n>1$. Let $\Gamma_{\bB_0}$ be the group
generated by inversions through each $(n+1)$-sphere $\Sigma^{n+1}_{j}=\partial B^{n+2}_{j}$,  whose Kleinian limit set is $\Lambda(\Gamma_{\bB_0})$,  where $B^{n+2}_{j}\in \bB_0$ $(j=1,2,\ldots,k)$.  Consider the inverse limit space $\mathcal{K}$ that is an $n$-dimensional wild knot invariant by a Schottky group. The purpose of this section is to construct cyclic branched covers of $\mathbb{S}^{n+2}$ along $\mathcal{K}$.\\

\noindent Since $\mathcal{K}$ is locally contractible, it follows by Alexander duality \cite{Hatcher} that
\[H^1(\bS^{n+2}{\smallsetminus}\mathcal{K},\Z)\cong H_{n}(\mathcal{K},\Z)\cong \Z.\]
Since $H^1(\bS^{n+2}{\smallsetminus}\mathcal{K},\Z)$ is isomorphic to the group of homotopy classes of maps of $\bS^{n+2}{\smallsetminus}\mathcal{K}$ into the circle $\bS^1$, it follows that $H^1(\bS^{n+2}{\smallsetminus} \mathcal{K},\Z)$ is generated by a non homotopically trivial map 
$h:\bS^{n+2}{\smallsetminus}\mathcal{K}\to\bS^1$. As a consequence, we obtain an infinite cyclic covering  $\tilde{\mathcal{K}}$ of $\mathbb{S}^{n+2}{\smallsetminus}\mathcal{K}$ such that the group $\mathbb{Z}$ acts freely as the group of covering transformations of $\tilde{\mathcal{K}}$ and the quotient space $\tilde{\mathcal{K}}/\mathbb{Z}$ can be identified with  
$\mathbb{S}^{n+2}{\smallsetminus}\mathcal{K}$. Furthermore, we can assume, that $h$ is a smooth map and then $h^{-1}(\theta)$, for $\theta\in\mathbb{S}^1$ a regular value of $h$, is a Seifert surface.

\begin{main3}\label{main3}
Let $B(K,\bB_0)$ be a $k$-beaded necklace subordinated
to the non-trivial, smooth $n$-knot $K$, and let  $\mathcal{K}$ be the corresponding wild $n$-knot. Then, for each integer $q$, there exists a $q$-fold cyclic branched cover 
$\Psi:\mathbb{M}_q\rightarrow \mathbb{S}^{n+2}$ along $\mathcal{K}$ such that $\mathbb{M}_q$ is  a compact and connected space.
\end{main3}

\noindent {\it Proof.} Suppose that $\bB_0=\bigcup_{j=1}^k B^{n+2}_j$, where $B^{n+2}_j$ is a closed, Euclidean $(n+2)$-ball, for $j=1,2,\ldots k$ and let $\Gamma_{\bB_0}$ be the group generated by inversions $I_j$ through each $(n+1)$-sphere $\Sigma^{n+1}_{j}=\partial B^{n+2}_{j}$ such that its Kleinian limit set is $\Lambda(\Gamma_{\bB_0})$. Consider the corresponding infinite cyclic cover 
$\tilde{P}:\tilde{\mathcal{K}}\rightarrow (\mathbb{S}^{n+2}{\smallsetminus}\mathcal{K})$ with fiber the Seifert 
surface $\tilde{S}$.  Thus, we can think $\tilde{\mathcal{K}}$ as the union of copies of the space 
$\mathbb{S}^{n+2}{\smallsetminus}\mathcal{K}$ cut open along $\tilde{S}$ (called $H$) and identifying these copies of the cut open space in a suitable way (see \cite{rolfsen}, pp.129--130).\\

\noindent Consider the subgroup $q\mathbb{Z}$ of $\mathbb{Z}$. Then there exists a covering map 
$\tilde{P}_q:\tilde{\mathcal{K}}_q\rightarrow (\mathbb{S}^{n+2}{\smallsetminus}\mathcal{K})$ where $\tilde{\mathcal{K}}_q=\tilde{\mathcal{K}}/\mathbb{Z}_q$. Now, let 
$R=K{\smallsetminus}\bigcup_{j=1}^k (K\cap B_j)$, so if $R_m=R_{m-1}\bigcup I_{j_m}I_{j_{m-1}} \cdots I_{j_1}(R)$, then the tame subset 
$\mathcal{K}{\smallsetminus}\Lambda(\Gamma_{\bB_0})$ is the direct limit space 
$\{R_m,\,|\,R_m\hookrightarrow R_{m+1}\}$, \emph{i.e.}, 
$\mathcal{K}{\smallsetminus}\Lambda(\Gamma_{\bB_0})=\bigcup_{m} R_m=\mathcal{R}$. Observe that 
$\mathcal{K}{\smallsetminus}\Lambda(\Gamma_{\bB_0})$ is a tame subset of $\mathbb{S}^{n+2}$, then we can extend this covering map to a $q$-fold cyclic cover 
$\psi_q:{\mathcal{M}}_q\rightarrow \mathbb{S}^{n+2}{\smallsetminus}\mathcal{K}$ branched over the tame set $\mathcal{R}=\mathcal{K}{\smallsetminus}\Lambda(\Gamma_{\bB_0})$, where 
$\mathcal{M}_q$ is obtained by pasting $q$ copies of $H\cup \mathcal{R}$ cyclically around $\mathcal{R}$.\\

\noindent Notice that ${\mathcal{M}}_q$ is a connected $(n+2)$-manifold with uncountable many ends. On the other hand, $\Lambda(\Gamma_{\bB_0})$ is a compact, totally disconnected subset of $\mathbb{S}^{n+2}$. Consider the inclusion map  
$j:(\mathbb{S}^{n+2}{\smallsetminus}\Lambda(\Gamma_{\bB_0}))\rightarrow \mathbb{S}^{n+2}$. Then the composition  $j\circ \psi$ is a spread, so there exists a unique Fox completion $\Psi:\mathbb{M}_q\rightarrow \mathbb{S}^{n+2}$ of the spread $j\circ \psi$, for more details see \cite{fox}, \cite{M}. This completion consists of compactifying each end, hence $\mathbb{M}_q={\mathcal{M}}_q\cup \Lambda(\Gamma_{\bB_0})$. The map $\Psi$ is induced by a cyclic action by homeomorphism in the topological space $\mathbb{M}_q$, hence the branched cover $\Psi:\mathbb{M}_q\rightarrow \mathbb{S}^{n+2}$ \emph{is a $q$-fold cyclic cover branched over the knot} $\mathcal{K}$. $\square$

\begin{rem}
As mentioned at the end of the introduction, in general, it is a complicated problem to prove that the Freudenthal compactification $\mathbb{M}$ of the $(n+2)$-manifold  ${\mathcal{M}}$ is an $(n+2)$-manifold
(for more details see \cite{fox} and \cite{M}). It is a challenging and beautiful problem to give conditions
on a wild knot in $\bS^{n+2}$, so that a branched cyclic covering over it is a topological manifold to obtain examples as those of Montesinos (\cite{M}).  
\end{rem}

\end{document}